\documentclass[10pt]{amsart}

\usepackage{eucal}
\usepackage{amscd}
\usepackage[all,cmtip]{xy}
\usepackage{rotating}
\usepackage{hyperref,mathrsfs}
\usepackage{tikz}

\usepackage{amsmath,amssymb}

\newtheorem{theorem}{Theorem }[section]

\newtheorem{lemma}[theorem]{Lemma}

\newtheorem{corollary}[theorem]{Corollary}

\theoremstyle{definition}
\newtheorem{definition}[theorem]{Definition}
\newtheorem{remark}[theorem]{Remark}

\def\I{\mathbf{I}}

\def\eop{\hspace*{\fill}{\footnotesize$\blacksquare$}}
\newcommand{\Aut}{\mathrm{Aut}}
\newcommand{\id}{\mathbf{1}}

\newcommand{\wis}[1]{{\text{\em \usefont{OT1}{cmtt}{m}{n} #1}}}

\newcommand{\mO}{\mathcal{O}}
\newcommand{\A}{\mathbb{A}}
\newcommand{\B}{\overline{\mathbf{B}}}

\newcommand{\hP}{\mathbb{P}}
\newcommand{\mD}{\mathcal{D}}
\newcommand{\Fun}{\mathbb{F}_1}
\newcommand{\Spec}{\wis{Spec}}
\newcommand{\Z}{\mathbb{Z}}
\newcommand{\Proj}{\wis{Proj}}
\newcommand{\mE}{\mathcal{E}}

\def\doubleprod#1#2{\ooalign{$#1\prod$\cr$#1\coprod$\cr}}
\DeclareMathOperator*{\Rprod}{\mathpalette\doubleprod\relax}

\newcommand{\fp}{\frak{p}}
\newcommand{\rd}{\mathrm{d}}

\newcommand{\mC}{\mathcal{C}}

\newcommand{\mF}{\mathcal{F}}

\newcommand{\fP}{\mathbf{P}}

\newcommand{\mS}{\mathcal{S}}

\newcommand{\F}{\mathbb{F}}
\newcommand{\mA}{\mathcal{A}}
\newcommand{\C}{\mathbb{C}}

\newcommand{\prf}{\textit{Proof. }}
\newcommand{\Sch}{\wis{Sch}}
\newcommand{\bL}{\mathbb{L}}
\newcommand{\uL}{\underline{\mathbb{L}}}
\usetikzlibrary{matrix}
\usetikzlibrary{arrows}
\usepackage{pgf}
\usepackage{lscape}
\usepackage{listings}

\usepackage{enumitem}

\title[Deitmar schemes and zeta functions]{Deitmar schemes, graphs and zeta functions}


\author{Manuel M\'{e}rida-Angulo and Koen Thas}

\address{{Ghent University},
{Department of Mathematics},
{Krijgslaan 281, S25, B-9000 Ghent, Belgium,}
\texttt{manmerang@gmail.com}
\texttt{koen.thas@gmail.com}}

\date{}

\begin{document}
\maketitle

\begin{abstract}
In \cite{KT-Japan} it was explained how one can naturally associate a Deitmar scheme (which is a scheme defined over the field with one element, $\Fun$) to a so-called ``loose graph'' (which is a generalization of a graph). Several properties of the Deitmar scheme can be proven easily from the 
combinatorics of the (loose) graph, and known realizations of objects over $\Fun$ such as combinatorial $\Fun$-projective and $\Fun$-affine spaces
exactly depict the loose graph which corresponds to the associated Deitmar scheme. In this paper, we first modify the construction of {\em loc. cit.}, and show that Deitmar schemes
which are defined by finite trees (with possible end points) are ``defined over $\Fun$'' in Kurokawa's sense; we then derive a precise formula
for the Kurokawa zeta function for such schemes (and so also for the counting polynomial of all associated $\mathbb{F}_q$-schemes). 
As a corollary, we find a zeta function for all such trees which contains information such as the number of inner points and the spectrum of degrees, and  which is thus very different than Ihara's zeta function (which is trivial in this case).  Using a process called ``surgery,'' we show that one can determine the zeta function of a general loose graph and its associated $\{$ Deitmar, Grothendieck $\}$-schemes  in a number of steps, eventually reducing the calculation essentially to trees. 
We study a number of classes of examples of loose graphs, and introduce the {\em Grothendieck ring of $\Fun$-schemes} along the way in order to perform the calculations.
Finally, we compare the new zeta function to Ihara's zeta function for graphs in a number of examples, and include a computer program for performing more tedious calculations.
\end{abstract}

\setcounter{tocdepth}{1}
\tableofcontents

\medskip
\section{Introduction}

\subsection{Deninger's program}

Let $\mC$ be a nonsingular absolutely irreducible algebraic curve over the finite field $\F_q$; its arithmetic zeta function is
\begin{equation}
\zeta_{\mC}(s) = \prod_{\fp}\frac{1}{1 - N(\fp)^{-s}},
\end{equation}
where $\fp$ runs through the closed points of $\mC$ and $N(\cdot)$ is the norm map. Fix an algebraic closure $\overline{\F_q}$ of $\F_q$ and let $m \ne 0$ be a positive integer; we have the following Lefschetz formula for the number $\vert \mC(\F_{q^m}) \vert$ of rational points over $\F_{q^m}$:

\begin{equation}
\vert \mC(\F_{q^m}) \vert = \sum_{\omega = 0}^2(-1)^{\omega}\mathrm{Tr}(\mathrm{Fr}^m \Big| H^\omega(\mC)) = 1 - \sum_{j = 0}^{2g}\lambda_j^m + q^f,
\end{equation}
where $\mathrm{Fr}$ is the Frobenius endomorphism acting on the \'{e}tale $\ell$-adic cohomology of $\mC$, the $\lambda_j$s are the eigenvalues of this action, and $g$ is the genus of the curve. It is not hard to show that we then have a ``weight decomposition''

\begin{equation}
\zeta_{\mC}(s) = \prod_{\omega = 0}^2\zeta_{h^{\omega}(\mC)}(s)^{(-1)^{\omega - 1}} 
		= \frac{\prod_{j = 1}^{2g}(1 - \lambda_jq^{-s})}{(1 - q^{-s})(1 - q^{1 - s})} \nonumber \\
\end{equation}

\begin{eqnarray} 	
	 = \frac{\mbox{\textsc{Det}}\Bigl((s\cdot\id - q^{-s}\cdot\mathrm{Fr})\Bigl| H^1(\mC)\Bigr)}{\mbox{\textsc{Det}}\Bigl((s\cdot\id - q^{-1}\cdot\mathrm{Fr})\Bigl| H^0(\mC)\Bigr.\Bigr)\mbox{\textsc{Det}}\Bigl((s\cdot\id - q^{-s}\cdot\mathrm{Fr})\Bigl| H^2(\mC)\Bigr.\Bigr)}.
	\end{eqnarray}
Here the $\omega$-weight component is the zeta function of the pure weight $\omega$ motive $h^{\omega}(\mC)$ of $\mC$.

Recalling the analogy between integers and polynomials in one variable over finite fields, 
Deninger | in a series of works \cite{Deninger1991,Deninger1992,Deninger1994} |  gave a description of conditions on a conjectural category of motives that would admit a translation of Weil's proof of the Riemann Hypothesis for function fields of projective curves over finite fields $\F_q$ to the hypothetical curve $\overline{\Spec(\Z)}$. In particular, he showed that the following formula would hold:

\begin{equation}
\zeta_{\overline{\Spec(\Z)}}(s) = 2^{-1/2}\pi^{-s/2}\Gamma(\frac s2)\zeta(s)  
		=  \frac{\Rprod_\rho\frac{s - \rho}{2\pi}}{\frac{s}{2\pi}\frac{s - 1}{2\pi}} \overset{?}{=} \\
\end{equation}

	\begin{eqnarray} 	
	 \frac{\mbox{\textsc{Det}}\Bigl(\frac 1{2\pi}(s\cdot\id - \Theta)\Bigl| H^1(\overline{\Spec(\Z)},*_{\mathrm{abs}})\Bigr.\Bigr)}{\mbox{\textsc{Det}}\Bigl(\frac 1{2\pi}(s\cdot\id -\Theta)\Bigl| H^0(\overline{\Spec(\Z)},*_{\mathrm{abs}})\Bigr.\Bigr)\mbox{\textsc{Det}}\Bigl(\frac 1{2\pi}(s\cdot\id - \Theta)\Bigl| H^2(\overline{\Spec(\Z)},*_{\mathrm{abs}})\Bigr.\Bigr)}, \nonumber 
	\end{eqnarray}
where $\Rprod$ is the infinite {\em regularized product}, similarly
$\mbox{\textsc{Det}}$ denotes the {\em regularized determinant} (a determinant-like function of operators on infinite dimensional vector spaces), $\Theta$ is an ``absolute'' Frobenius endomorphism, and the $H^i(\overline{\Spec(\Z)},*_{\mathrm{abs}})$ are certain proposed cohomology groups. The $\rho$s run through the set of critical zeroes of the classical Riemann zeta. (In the formula displayed above, $\Spec(\Z)$ is compactified to $\overline{\Spec(\Z)}$ in order to see it as a projective curve.)

Combined with Kurokawa's work on multiple zeta functions (\cite{Kurokawa1992}) from 1992, the conjectural idea arose that there are motives $h^0$ (``the absolute point''), $h^1$ and $h^2$ (``the absolute Lefschetz motive'') with zeta functions
	\begin{equation}
	\label{eqzeta}
		\zeta_{h^w}(s) \ = \ \mbox{\textsc{Det}}\Bigl(\frac 1{2\pi}(s\cdot\id-\Theta)\Bigl| H^w(\overline{\Spec(\Z)},*_{\mathrm{abs}})\Bigr.\Bigr) 
	\end{equation}
for $w=0,1,2$. Deninger computed that 
\begin{equation}
\label{Denzeta}
\zeta_{h^0}(s)=s/2\pi\ \ \mbox{and}\ \  \zeta_{h^2}(s)=(s-1)/2\pi. 
\end{equation}

In his now famous lecture notes \cite{Manin}, Manin suggested to interpret  $h^0$ as $\Spec(\Fun)$, where 
$\Fun$ denotes the ``field with one element,'' and  $h^2$ as the affine line over  $\Fun$. Essentially from that moment, the search for a proof of the Riemann Hypothesis became a main motivation to look for a geometric theory over $\Fun$. 

And on the other hand, in this larger Grothendieck site of Deninger, the Cartesian product
\begin{equation}
\Spec(\Z) \times_{\Upsilon} \Spec(\Z) \times_{\Upsilon} \cdots \times_{\Upsilon} \Spec(\Z)
\end{equation}
over the deeper base $\Upsilon = \Fun$ eventually should make sense (as it is crucial in Weil's proof).\\

\subsection{$\Fun$-Schemes and their zeta functions}

By 2006, there was a scheme theory developed by Deitmar in  \cite{Deitmarschemes2} over $\Fun$ | related to Kato's ``log scheme theory'' \cite{Kato} | which is based on the observation that commutative rings over $\Fun$ could be imagined as commutative multiplicative  monoids (with an absorbing element). For these ``rings,'' one  can naturally define (prime) ideals, localization, a $\Spec$-construction, etc., and Deitmar has proven that a natural base extension
to $\Z$ of varieties in this context leads to toric varieties \cite{Deitmarschemes1}. (We refer to \S \ref{Deitsch} for more details and examples.) Other, often more general, scheme theories over $\Fun$ have been defined | see for instance \cite{Map,ChapOL} for an account | but in one way or the other, Deitmar schemes 
appear to be the core of such scheme theories. In any case Deitmar schemes are amongst  the most basic objects in $\Fun$-theory, and 
one needs to understand them as well as possible. (We refer to the monograph \cite{AA} for background on the object $\Fun$.)

There also is a promising theory for zeta functions of Deitmar schemes \cite{Deitmarschemes1} which agrees with Kurokawa's theory for zeta functions of schemes ``defined over $\Fun$'' \cite{Kurozeta}, and the latter objects are one of the main subjects of the present paper | see \S \ref{Kuro} for formal definitions. 
(After the Tate conjectures, they appear to be precisely those schemes which come with a mixed Tate motive, see \S \ref{Kuro}.)

In fact, impressively enough, Kurokawa's $\Fun$-zeta functions for $\A^n$ and $\hP^n$ | both spaces being defined over $\Fun$ | agree with Deninger's computation (\ref{Denzeta}) (ignoring the $\frac{1}{2\pi}$-factor).

\subsection{The combinatorial side}

Besides the Algebraic Geometry side of $\Fun$-theory, there is also the combinatorial-synthetic side: in an old paper \cite{anal} which pre-dates any of the papers we mentioned so far, Tits
described symmetric groups as ``Chevalley groups over $\Fun$,'' as limit objects of projective general linear groups over finite 
fields, and the corresponding geometric modules, which are just projective spaces over the same fields in this case, 
become complete graphs with the full subgraph structure.

In the note \cite{KT-Japan}, the second author united these two theories in the following way: starting with a ``loose graph'' $\Gamma$, which is a more general object than a graph, he associated a Deitmar scheme $\mS(\Gamma)$ with $\Gamma$ of which the closed points correspond to the vertices of the loose graph. Several fundamental properties of the Deitmar scheme can  be obtained easily from the 
combinatorics of the loose graph, such as connectedness, the automorphism group, etc. And very interestingly, it appears
that a number of combinatorial $\Fun$-objects (such as the combinatorial $\Fun$-projective space of above) are just loose graphs, and moreover,
the associated Deitmar schemes are precisely  the scheme versions in Deitmar's theory of these objects. Also, the automorphism groups of 
$\Gamma$ and $\mS(\Gamma)$ are isomorphic.

\subsection{The present paper}

In this paper, we merge the ideas of Kurokawa \cite{Kurokawa1992}, Deninger's calculation (\ref{Denzeta}) and Manin's interpretation of 
(\ref{Denzeta}) with the ideas of \cite{KT-Japan}, and
introduce a new zeta function for all (loose) graphs, and the $\{$ Deitmar/Grothendieck $\}$-schemes associated to them. 
The idea is, as in \cite{KT-Japan}, to associate a Deitmar scheme to a loose graph in a natural way (but not in the same way as in \cite{KT-Japan}), 
and to show that, after having applied Deitmar's ($\cdot \otimes_{\Fun}\mathbb{Z}$)-functor, the obtained Grothendieck schemes are defined over $\Fun$ in Kurokawa's sense. So they come with a Kurokawa zeta function. \\

We first slightly modify (generalize) the functor $\mS$ (and call it ``$\mF$''), so that it becomes more natural for our purposes.
The functor $\mF$ has to obey a tight set of rules in order to meet our purposes at the end. We mention a few of the important ones:
\begin{itemize}
\item[{\bf Rule} \#1]
The loose graphs of the affine and projective space Deitmar schemes (as described later on) should correspond to affine and projective space Deitmar schemes.
\item[{\bf Rule} \#2]
A vertex of degree $m$ should correspond locally to an affine space $\A^m$.
\item[{\bf Rule} \#3]
An edge without vertices should correspond to a multiplicative group.
\item[{\bf Rule} \#4]
``The loose graph is the map to gluing.''
\end{itemize}

Because of Rule \#1, the pictures of Tits and Kapranov-Smirnov of affine and projective spaces over $\Fun$ | which are loose graphs (see \cite{Chap1,Chap2} for details) |  are in agreement with our functor $\mF$.

We  show that Deitmar schemes
which are defined through $\mF$  by finite trees (with possible end points) give $\mathbb{Z}$-schemes which are ``defined over $\Fun$'' in Kurokawa's sense, and we derive a precise formula for the Kurokawa zeta function for such schemes (and so also for the counting polynomial of all associated $\mathbb{F}_q$-schemes) which can be read from the (loose) trees. 
As a corollary, we find a zeta function for all such loose trees which contains information such as the number of inner points and the spectrum of degrees, and  which is thus very different than Ihara's zeta function (which is trivial in this case). 

To calculate the Kurokawa zeta functions of schemes associated to general loose graphs, we develop a tedious procedure called ``surgery,'' which determines what happens to the zeta function of such a scheme, after one has replaced one particular edge (with two vertices $u$ and $v$), by two
edges, each  with one vertex ($u$, respectively $v$). Repeating this procedure a number of times, and keeping track of the differences in each step, one eventually ends up with a loose tree, and this is the case which was already settled. 

All calculations are done through the {\em Grothendieck ring of $\Fun$-schemes} (which is introduced for this purpose), and in order to 
develop surgery, we show that to know such a difference between zeta functions, one only has to understand what happens in a neighborhood of the edge ($uv$) which is ``resolved.'' 

We study surgery for a number of classes of particular loose graphs before passing to the general case, and explain  several subtleties and booby traps on examples. 

As calculating the zeta function for a loose graph is even very complicated for seemingly simple examples, 
we also include a computer program for this purpose.

\newpage
\section{The Ihara zeta function}
\label{Iharasec}

Let $\Gamma = (V,E)$ be a finite connected undirected graph with no vertices of degree $1$ (``end points''). The {\em rank} $r_{\Gamma}$ of $\Gamma$ is 
$\vert E \vert - \vert V \vert + 1$; it is the number of edges one has to delete from $\Gamma$ to obtain a spanning tree. Alternatively, one can also define $r_{\Gamma}$ as the rank of the fundamental group of $\Gamma$. Suppose that $r_{\Gamma} \geq 1$ | that is, $\Gamma$ is not a tree. 

Let $E = \{ e_1, \ldots, e_n\}$, and define a new oriented edge set of size $2\vert E\vert$ as follows (where below the edges of $E$ are arbitrarily oriented):
\begin{equation}
e_1,\ldots,e_n;\ \ \ e_{n + 1} = e^{-1}_1,\ldots, e_{2n} = e^{-1}_n. 
\end{equation}

Let $D = a_1a_2\cdots a_r$ be a directed path in $\Gamma$ (all the $a_i$ are edges and they are directed in the same direction). 
The {\em equivalence class} $[D]$ of $D$ is the set
\begin{equation}
[D] = \{ a_1a_2\cdots a_r, a_2a_3\cdots a_ra_1, \ldots, a_ra_1\cdots a_{r - 1}\}.
\end{equation}

The {\em length} $\nu(D)$ of $D$ is the number $r$.
We say $D$ is ``backtrackless'' if 
$a_{i + 1} \ne a_i^{-1}$ for all $i$, and it is ``tailless'' if $a_r \ne a_1^{-1}$. The path $D$ is {\em primitive} if $D \ne F^m$ for any postive integer $m \geq 2$ and any directed path $F$.

A {\em prime} $[P]$ for $\Gamma$ is an equivalence class of closed backtrackless tailless primitive (directed) paths in $\Gamma$.

The {\em Ihara zeta function} of $\Gamma$ is now defined as follows:
\begin{equation}
\zeta(u,\Gamma) := \prod_{[P]\ \ \mbox{prime}}(1 - u^{\nu(P)})^{-1}.
\end{equation}
Here $u \in \mathbb{C}$ with $\vert u \vert$ sufficiently small.

One way to calculate $\zeta(u,\Gamma)$ is through the Bass-Hashimoto formula
\begin{equation}
\zeta(u,\Gamma)^{-1} = (1 - u^2)^{r_{\Gamma} - 1}\mathrm{det}(\id - A_{\Gamma}u + Q_{\Gamma}u^2),
\end{equation}
where $\id$ is an identity matrix, $A_{\Gamma}$ the adjacency matrix of $\Gamma$, and $Q_{\Gamma}$ the diagonal matrix whose $j$-th
diagonal entry is ($-1 +$ degree of $j$-th vertex). (One has numbered the vertices in order to obtain $\id, A_{\Gamma}$ and $Q_{\Gamma}$; all are $(\vert V \vert \times \vert V \vert)$-matrices.) 

\begin{remark}{\rm
In general, $\zeta(u,\Gamma)^{-1}$ is a polynomial of degree $2\vert E\vert$.}
\end{remark}

Put $\vert E \vert = m$, and define a $(2m \times 2m)$-matrix $\mE$ (the ``edge adjacency matrix'') by letting the $ij$-th entry be $1$ if the terminal vertex of $e_i$ is the initial vertex of $e_j$, provided that $e_j \ne e_i^{-1}$. Otherwise, the entry is $0$. Then it can be shown that the Ihara zeta function of $\Gamma$ can also be calculated as
\begin{equation}
\zeta(u,\Gamma)^{-1} = \mathrm{det}(\id - u\mE). 
\end{equation}

In other words, the roots of $\zeta(u,\Gamma)^{-1}$ (with multiplicities) are the eigenvalues (with multiplicities) of $\mE$. So two graphs have the same Ihara zeta function{\em  if and only if they are isospectral with respect to the edge adjacency matrix}.\\

The following properties/values can be read from the Ihara zeta function of a (finite connected undirected) graph (with no vertices of degree $1$, and rank $\geq 1$):
\begin{itemize}
\item
whether it is bipartite or not;
\item
its number of vertices and edges;
\item
whether it is regular, and if so, its regularity degree and spectrum.
\end{itemize}

\newpage
\section{Schemes defined over $\Fun$, and their zeta functions}
\label{Kuro}

In \cite{Kurozeta}, Kurokawa says a scheme $X$ is of {\em $\Fun$-type}\index{scheme of $\Fun$-type} if its arithmetic zeta function $\zeta_X(s)$ can be expressed in the form
\begin{equation}
\zeta_X(s) = \prod_{k = 0}^n\zeta(s - k)^{a_k}
\end{equation}
with the $a_k$s in $\Z$. A very interesting result in \cite{Kurozeta} reads as follows:

\begin{theorem}
Let $X$ be a $\Z$-scheme. The following are equivalent.
\begin{itemize}
\item[{\rm (i)}]
\begin{equation}
\zeta_X(s) = \prod_{k = 0}^n\zeta(s - k)^{a_k}
\end{equation}
with the $a_k$s in $\Z$.
\item[{\rm (ii)}]
For all primes $p$ we have
\begin{equation}
\zeta_{X\vert \F_p}(s) = \prod_{k = 0}^n(1 - p^{k - s})^{-a_k}
\end{equation}
with the $a_k$s in $\Z$.
\item[{\rm (iii)}]
There exists a polynomial $P_X(Y) = \sum_{i = 0}^na_kY^k$ such that
\begin{equation}
\#X(\F_{p^m}) = P_X(p^m) 
\end{equation}
for all finite fields $\F_{p^m}$.
\end{itemize}
\end{theorem}

Kurokawa defines the {\em $\Fun$-zeta function}\index{$\Fun$-zeta function} of a $\Z$-scheme $X$ which is defined over $\Fun$ as 
\begin{equation}
\zeta_{X\vert \Fun}(s) :=  \prod_{k = 0}^n(s - k)^{-a_k}
\end{equation}
with the $a_k$s as above. Define, again as above, the {\em Euler characteristic}\index{Euler characteristic}
\begin{equation}
\#X(\Fun) := \sum_{k = 0}^na_k.
\end{equation}

The connection between $\Fun$-zeta functions and arithmetic zeta functions is explained in the following theorem, taken from \cite{Kurozeta}.

\begin{theorem}
Let $X$ be a $\Z$-scheme which is defined over $\Fun$. Then
\begin{equation}
\zeta_{X\vert \Fun}(s) =  \lim_{p \longrightarrow 1}\zeta_{X\vert \F_p}(s)(p - 1)^{\# X(\Fun)}.
\end{equation}
Here, $p$ is seen as a complex variable (so that the left hand term is the leading coefficient of the Laurent expansion of $\zeta_{X \vert \Fun}(s)$ around $p = 1$).
\end{theorem}

For affine and projective spaces, we obtain the following zeta functions (over $\Z$, $\F_p$ and $\Fun$, with $n \in \mathbb{N}^\times$):
\begin{eqnarray}
\zeta_{\A^n\vert \Z}(s) &= &\zeta(s - n);\nonumber \\
\zeta_{\A^n\vert \F_p}(s) &= &\frac{1}{1 - p^{n - s}};\nonumber \\
\zeta_{\A^n\vert \Fun}(s) &= &\frac{1}{s - n},
\end{eqnarray}
and
\begin{eqnarray}
\zeta_{\fP^n\vert \Z}(s) &= &\zeta(s)\zeta(s - 1)\cdots\zeta(s - n);\nonumber \\
\zeta_{\fP^n\vert \F_p}(s) &= &\frac{1}{(1 - p^{-s})(1 - p^{1 - s})\cdots(1 - p^{n - s})};\nonumber \\
\zeta_{\fP^n\vert \Fun}(s) &= &\frac{1}{s(s - 1)\cdots(s - n)}.\\
\end{eqnarray}

\medskip
\section{Deitmar schemes}
\label{Deitsch}

Define an {\em $\Fun$-ring} to be a commutative monoid with an absorbing element $0$.

For the definition of Deitmar scheme, we refer to \cite{Deitmarschemes2} or \cite{Chap2}. Let us just mention
that one defines an {\em affine Deitmar scheme} $\Spec(A)$ similarly as an affine Grothendieck scheme (by building a Zariski-type
topology  and structure sheaf of $\Fun$-rings on the set of monoidal prime ideals of the commutative monoid with zero $A$). A general {\em Deitmar scheme} then is  a {\em monoidal space} (a topological space endowed with a sheaf of $\Fun$-rings), locally isomorphic to affine Deitmar schemes. (Deitmar schemes  are sometimes also called {\em $\mD_0$-schemes} or {\em $\mathcal{M}_0$-schemes} in some papers.)

We list some of the most important examples below.

\medskip
\subsection{Polynomial rings}

Define
\begin{equation}
\mathbb{F}_1[X_1,\ldots,X_n] := \{0\} \cup  \{X_1^{u_1}\ldots X_n^{u_n} \vert u_j \in \mathbb{N}\}, 
\end{equation}\index{$\mathbb{F}_1[X_1,\ldots,X_n]$}
that is, the union of $\{0\}$ and the (abelian) monoid generated by the $X_j$.

\medskip
\subsection{Affine space}

Let $A=\Fun[X_1,\ldots,X_n]$.  Denote $\Spec({\Fun[X_1,\ldots,X_n]})$ by $\A_{\Fun}^n$ and call it the \emph{$n$-dimensional affine space over $\Fun$}\index{affine!space over $\Fun$}. The $\ne (0)$ prime ideals of $A$ are of the form 
$\fp_I=\bigcup_{i\in I}(X_i)$, where $I$ is a subset of $\{1,\dotsc,n\}$ and $(X_i) = X_iA=\{X_ia\mid a\in A\}$.

\medskip
\subsection{$\Proj$-schemes}

In \cite{NotesI} the second author introduced the $\Proj$-scheme construction for Deitmar schemes (see also \cite{Chap2}). We recall the procedure in a nutshell.

 \subsubsection{Monoid quotients}
 
 Let $M$ be a commutative unital monoid (with $0$), and $I$ an ideal of $M$. We define the monoidal quotient $M/I$ to be the set $\{ [m] \in M \vert m \in M \}/([m] = [0]\ \mbox{if}\ m \in I)$.
 (When $R$ is a commutative ring and $J$ an ideal, then the ring quotient $R/J$ induces the monoidal quotient on $R, \times$.)
 
  \subsubsection{The $\Proj$-construction}
 
 Consider the $\mathbb{F}_1$-ring $\mathbb{F}_1[X_0,X_1,\ldots,X_m]$, where $m \in \mathbb{N}$. Since any polynomial is 
 homogeneous in this ring, we have a natural grading
 
 \begin{equation}
 \mathbb{F}_1[X_0,\ldots,X_m] = \bigoplus_{i \geq 0}R_i = \coprod_{i \geq 0}R_i,
 \end{equation}
 where $R_i$ consists of the elements of $\mathbb{F}_1[X_0,X_1,\ldots,X_m]$ of total degree $i$, for $i \in \mathbb{N}$. 
 The {\em irrelevant ideal}\index{irrelevant ideal}\index{\mathrm{Irr}} is defined as
 
 \begin{equation}
\mathrm{Irr} = \{0\} \cup \coprod_{i \geq 1}R_i.
 \end{equation}
 
 (It is just the monoid minus the element $1$.) 
 Now $\Proj(\mathbb{F}_1[X_0,\ldots,X_m]) =: \Proj(\mathbb{F}_1[\mathbf{X}])$\index{$\Proj(\mathbb{F}_1[\mathbf{X}])$} consists, as a set, of the prime ideals of  $\mathbb{F}_1[X_0,X_1,\ldots,X_m]$
 which do not contain $\mathrm{Irr}$ (so only $\mathrm{Irr}$ is left out of the complete set of prime ideals). 
 The closed sets of the (Zariski) topology on this set are defined as usual: for any ideal $I$ of  $\mathbb{F}_1[X_0,X_1,\ldots,X_m]$, we define
 
 \begin{equation}
 V(I) := \{ \fp \vert \fp\in \Proj(\mathbb{F}_1[\mathbf{X}]),\ \ I \subseteq \fp \},
 \end{equation}
 where $V(I) = \emptyset$ if $I = \mathrm{Irr}$ and $V(\{0\}) = \Proj(\mathbb{F}_1(\mathbf{X}))$,
 the open sets then being of the form
  \begin{equation}
 D(I) := \{ \fp \vert \fp\in \Proj(\mathbb{F}_1[\mathbf{X})],\ \ I \not\subseteq \fp \}.
 \end{equation}
 
 It is obvious that $\Proj(\mathbb{F}_1[\mathbf{X}])$ is a Deitmar scheme.
  Each ideal $(X_i)$ defines an open set $D((X_i))$ such that the restriction of the scheme to this set is isomorphic to $\Spec(\mathbb{F}_1[\mathbf{X}_{(i)}])$, where $\mathbf{X}_{(i)}$ is $X_0,X_1,\ldots,X_m$ with $X_i$ left out.\\

 More generally, suppose $M$ is any commutative unital monoid (with $0$) with a  grading
 \begin{equation}
M = \coprod_{i \geq 0}M_i,
 \end{equation}
 where the $M_i$ are the sets with elements of total degree $i$ (for $i \in \mathbb{N}$), and let, as above, the
{\em irrelevant ideal}\index{irrelevant ideal} be $\mathrm{Irr} = \{0\} \cup \coprod_{i \geq 1}M_i$\index{$\mathrm{Irr}$}. Define the topology $\Proj(M)$\index{$\Proj(M)$} as before (noting that homogeneous (prime) ideals are the same as ordinary monoidal (prime) ideals here). For an open 
$U$, define $\mO_M(U)$ as consisting of all functions
 \begin{equation}
f: U \longrightarrow \coprod_{\fp \in U}M_{(\fp)},
\end{equation}
where $M_{(\fp)}$ is the subset of $M_{\fp}$ of fractions of elements with the same degree, 
for which $f(\fp) \in M_{(\fp)}$ for each $\fp \in U$, and such that there exists a neighborhood $V$ of $\fp$ in $U$, and  elements $u, v \in M$, for which
$v \not\in \mathfrak{q}$ for every $\mathfrak{q} \in V$, and $f(\mathfrak{q}) = \frac{u}{v}$ in $M_{(\mathfrak{q})}$.

In this way we obtain a sheaf of $\mathbb{F}_1$-rings on $\Proj(M)$ making it a Deitmar scheme.

  \begin{remark}[Closed points]
 \label{combsch}
 {\rm
 Note that the closed points and projective sublines of $\Proj(\mathbb{F}_1[X_0,\ldots,X_m])$ form a complete graph on $m + 1$ vertices, so we 
 can easily switch between combinatorial $\mathbb{F}_1$-projective spaces and
 $ \Proj(\mathbb{F}_1[\mathbf{X}])$-schemes.\\
 }
\end{remark}

\medskip
\section{Loose graphs and the functor $\mS$}

\subsection{Loose graphs}

Define a {\em loose graph}\index{loose!graph} to be a triple $(V,E,\I)$, where as in graph theory $V$ is a set of {\em vertices}, $E$ is a set of {\em edges} ($V \cap E = \emptyset$), and $\I$
is a symmetric relation on $(V \times E) \cup (E \times V)$ which indicates when a vertex and an edge are incident,
with the additional property that {\em each edge is incident with at most two distinct vertices}. 
It relaxes the definition of graphs, in that an edge can now also have one, or even no, point(s). Also, we do not allow loops, and the geometry is undirected. Usually we will consider connected loose graphs.

 \medskip
\subsection{Embedding theorem} 
\label{EmbThm}
 
 Let $\Gamma$ be a loose graph. The embedding theorem of \cite{KT-Japan} observes that $\Gamma$ can be seen as a 
 subgeometry of the combinatorial projective $\Fun$-space $\mathbf{P}(\Gamma)$, called the {\em ambient} space, by simply 
 adding a second vertex on each edge which contains only one vertex, and then constructing the complete graph on the 
 total set of vertices. We will use the same notation $\mathbf{P}(\Gamma)$ for the associated projective space scheme
 (and the scheme over which the space is considered will be specified accordingly).

 \medskip
\subsection{The functor $\mS$} 
 
 Now let $\Gamma = (V,E,\I)$ be  a not necessarily finite graph. We will give a ``patching" argument as follows.

Consider $\fP = \fP(\Gamma)$, and note that since $\Gamma$ is a graph, $\fP \setminus \Gamma$ | when $\fP$ is considered as a graph | is just a 
set $S$ of edges. Let $\mu$ be arbitrary in $S$, and let $z$ be one of the two (closed) points on $\mu$ in $\fP = \Proj(\mathbb{F}_1[X_i]_{i \in V})$  (recall Remark \ref{combsch}).
Suppose that in the projective space $\fP$, $z$ is 
defined by the ideal generated by the polynomials
\begin{equation}
X_i,\ \ i \in V, i \ne j = j(z).
\end{equation}

Let $\fP(z)$ be the complement in $\fP$ of $z$; it is a hyperplane defined by $X_j = 0$ (and it forms a complete graph on all the points but $z$).
Denote the corresponding closed subset of $\Proj(\mathbb{F}_1[X_i]_{i \in V})$ by $C(z)$.
Let $z' \ne z$ be the other point of the edge $\mu$ corresponding to the index $j' = j(z') \in V$. Define the subset $\fP(z') = \fP \setminus \{z'\}$ of $V$, and denote the corresponding closed subset by $C(z')$. 
Finally, define
\begin{equation}
C({\mu}) = C(z) \cup C(z').
\end{equation}

It is also closed in $\Proj(\mathbb{F}_1[X_i]_{i \in V})$, and the corresponding closed subscheme is the projective space $\fP$ ``without the edge $\mu$''; the coordinate ring is $\mathbb{F}_1{[X_i]}_{i \in V}/I_{\mu}$ (where $(X_jX_l) =: I_{\mu}$) and its scheme is the $\Proj$-scheme defined by this ring.
Now introduce the closed subset
\begin{equation}
C(\Gamma) =  \cap_{\mu \in S}C(\mu).
\end{equation}

Then $C(\Gamma)$ defines a closed subscheme $S(\Gamma)$ which corresponds to the graph $\Gamma$.  We have 
\begin{equation}
S(\Gamma) = \Proj(\mathbb{F}_1{[X_i]}_{i \in V}/\cup_{\mu \in S}I_{\mu}).
\end{equation}

\begin{remark}[Edges and relations]{\rm
In this presentation, an edge corresponds to a relation, and we construct a coordinate ring for $\Theta(\Gamma) = S(\Gamma)$ by deleting all relations of the ambient space $\fP(\Gamma)$ which  are defined by edges in the complement of $\Gamma$.}
\end{remark}

A similar construction can be done for loose graphs, \cite{NotesI}.

\medskip
\subsection{Some properties}

In \cite{KT-Japan}, it is shown that certain properties of Deitmar schemes arising from loose graphs can be easily verified on the loose graphs. We mention some result.

 The next theorem shows that the automorphism group of projective spaces from the incidence geometrical point of view, which we denote by $\Aut_{\mathrm{synth}}(.)$, coincides with the automorphism group from the point of view of $\mathbb{F}_1$-schemes, denoted $\Aut_{\mathrm{sch}}(.)$.

\begin{theorem}[\cite{KT-Japan}]
Let $\fP$ be a projective space over $\mathbb{F}_1$, and let $\Proj(\mathbb{F}_1[X_i]_{i \in \fP})$ be the corresponding projective scheme.
Then we have
\begin{equation}
\Aut_{\mathrm{synth}}(\fP) \cong \Aut_{\mathrm{sch}}(\Proj(\mathbb{F}_1[X_i]_{i \in \fP})).
\end{equation}
\end{theorem}

A similar proof (considering the action on the ideals that correspond to the ``directions'' instead of the closed points) leads to the same theorem for affine spaces.
 
\begin{corollary}[\cite{KT-Japan}]
Each group $H$ is the full automorphism group of some Deitmar scheme.
\end{corollary} 
 
 Denote the category of loose (undirected, loopless)
graphs and natural morphisms by $\mathbf{LG}$\index{$\mathbb{LG}$}.

\begin{theorem}[\cite{KT-Japan}]
For any element $\Gamma \in \mathbf{LG}$, we have that
 \begin{equation}
 \Aut(\Gamma)_{\mathrm{synth}} \cong \Aut(S(\Gamma))_{\mathrm{sch}}.
 \end{equation}
 \end{theorem}
 
\begin{theorem}[\cite{KT-Japan}]
A  loose scheme $S(\Gamma)$ is connected if and only if the  loose graph $\Gamma$ is connected.
\end{theorem}

\newpage
\section{Modifying the functor $\mS$}
\label{ppme}

Let $B$ be the complete graph on three vertices minus one edge: by the embedding theorem, it embeds in projective
$2$-space over $\Fun$. Call the vertex of degree two $x$, and call the others $y$ and $z$. It contains the graph $A$ corresponding to $\Spec(\Fun[X_1,X_2])$ as an  open set. On the level of $\mD$-schemes, what we obviously want is $\mS(B)$ to be the projective plane over $\Fun$ minus 
the multiplicative group (since we delete a projective line minus two points). And very importantly, we also want that 
if $\Gamma \subset \widetilde{\Gamma}$ is a strict inclusion of loose graphs, $\mS(\Gamma)$ also is a proper subscheme
of $\mS(\widetilde{\Gamma})$. In fact, we will introduce this property as one of our axioms to define the modified functor (which we still call $\mS$ for now).

\medskip
\begin{itemize}
\item[COV]
If $\Gamma \subset \widetilde{\Gamma}$ is a strict inclusion of loose graphs, $\mS(\Gamma)$ also is a proper subscheme
of $\mS(\widetilde{\Gamma})$.
\item[LOC-DIM]
If $x$ is a vertex of degree $m \in \mathbb{N}^\times$ in $\Gamma$, then there is a neighborhood $\Omega$ of $x$ in 
$\mS(\Gamma)$ such that $\mS(\Gamma)_{\vert \Omega}$ is an affine space of dimension $m$.
\end{itemize}

\medskip
Now consider an open cover $A_1,A_2,A_3$ of $\mS(B)$, where $A_1$ contains $y$ and is defined by taking out the closed subscheme defined by the projective line $xz$, $A_2$ contains $x$ and is defined by taking out the closed points $y$ and $z$, and $A_3$ contains $z$ and is similarly defined as $A_1$. With the same indices, we obtain affine schemes $\Spec(E_i)$ which in case of $i = 1,3$ is the absolute line and in case $i = 2$ the affine plane over $\Fun$.  Applying Deitmar's  base-extension to $\mathbb{Z}$, $\mS(B) \otimes_{\Fun} \mathbb{Z}$ is defined 
by lifting the $\Spec(E_i)$ to $\mathbb{Z}$, the gluing being governed by what happens on the $\mD$-level (see also \S \ref{glue}).

For any finite field $k = \mathbb{F}_q$, we obtain a  scheme $\mS(B)(k) \hookrightarrow \fP^2(\mathbb{F}_q)$ which is covered by three opens, in which two 
affine lines and one affine plane are induced, and the gluings are defined by the drawing (cf. \S \ref{glue}), that is, $\fP^2(\mathbb{F}_q)$ without a projective line minus two points.  Choosing homgeneous coordinates $X, Y, Z$ in $\fP^2(\mathbb{F}_q)$ such that the aforementioned projective line is defined by $Z  = 0$, we obtain the quasi-projective variety $V(k)$ of which the complement is defined by the property
\begin{equation}
Z = 0\ \ \Longrightarrow\ \ XY \ne 0.
\end{equation}

Whence $V(k)$ is given by the constructible set
\begin{equation}
(Z \ne 0)\ \ \vee\ \ (XY = 0),
\end{equation}
that is, the union of the $X$- and $Y$-coordinate axes and the affine plane with as line at infinity $Z = 0$.

\begin{remark}
\label{conf}
{\rm
Although, going from $A$ to $B$, we add the points on the edges through $x$, this does not imply, as we have seen, that once we have lifted 
the $\mD$-schemes $\mS(A)$ and $\mS(B)$ to $\mathbb{Z}$-schemes or $\mathbb{F}_q$-schemes, all the points of the lines in the 
affine plane corresponding to $A$ get an extra point. In fact, an affine plane over $\Fun$ is just an absolute point with two directions, 
so these (two) lines correspond to canonical coordinate axes. Once lifted to $\mathbb{Z}$ or $\mathbb{F}_q$, only the (two) lines corresponding to these axes get an extra point. 
}
\end{remark}

As in \cite{KT-Japan}, we also want the following two properties to hold:
\medskip
\begin{itemize}
\item[C]
any vertex of a loose graph $\Gamma$ defines a closed point of $\mS(\Gamma)$, and vice versa;
\item[PL]
any edge of $\Gamma$ which contains two distinct vertices corresponds to a closed sub projective line of $\mS(\Gamma)$.
\end{itemize}

Both are in fact instances of the following more general property: 
\medskip
\begin{itemize}
\item[CO]
If $K_m$ is a sub complete graph on $m$ vertices in $\Gamma$, then $\mS(K_m)$ is a closed sub projective space 
of dimension $m - 1$ in $\mS(\Gamma)$. 
\end{itemize}

\medskip
\subsection{Open and closed sets | Caution!}

Let $\Gamma$ be any finite loose graph, and let $L$ and $L'$ be two different sub complete graphs on two vertices. After passing to 
$\mS(\Gamma)$, we obtain different closed projective lines $\mS(L)$ and $\mS(L')$. Unfortunately, the subgraph of $\Gamma$ defined by the union 
$L \cup L'$ does not necessarily define a closed subset of $\mS(\Gamma)$ (although it looks rather closed in the loose graph); in the example in the beginning of this section, for example, that union defined a proper constructible set. 

The reason of this possible confusion, as already mentioned in Remark \ref{conf}, is that although, for instance in the aforementioned example, 
$L \cup L'$ looks like the union of two projective lines meeting in a point, it is in fact a projective plane minus a multiplicative group. Unfortunately (again), we only see this after the scheme has acquired enough flesh | that is, after base change to ``real'' fields, cf. Remark \ref{conf}. It is essentially a corollary of Axiom (COV): the union $L \cup L'$ (with $L, L'$ still meeting in a point) contains the loose graph of the $\Fun$-affine plane, 
so must define something $2$-dimensional. 

On the other hand, if $L$ and $L'$ would not meet in some loose graph $\Gamma$, then the union {\em does} define a closed set (by simply multiplying equations). In general, the same holds for a finite number of sub projective lines in general position.

Similarly, one has to be careful with finite unions of sub complete grahs. But finite unions of vertices always 
define closed sets | call this property ``FUCP'' for further reference.

\begin{remark}{\rm
If (COV) is not asked, a given loose graph scheme could (more or less) naturally ascend to {\em several} schemes over a fixed finite field $\mathbb{F}_p$.  
We will come back on this matter in a later paper.
}
\end{remark}

\medskip
\subsection{Lifting of general ``loose stars''}

After the analysis on the example in this section, and taken the main axioms into account, it is now easy to write down
the constructible set for the $\mathbb{Z}$-scheme corresponding to a general ``loose star graph.'' Let $\Gamma$ have a vertex $v$
which is incident with $m \geq 2$ edges, of which $\ell \leq m$, $\ell \geq 0$, have a second vertex (and let $\Gamma$ contain no other vertices and edges).
Then after choosing homogeneous coordinates in $\mathbf{P}^{m - 1}(\mathbb{Z})$, the corresponding constructible 
set in $\mathbb{Z}[X,X_1,\ldots,X_m]$ is 

\begin{equation}
(X \ne 0)\ \ \vee\ \ (\{ X_iX_j = 0 \vert i \ne j; 1 \leq i, j \leq \ell \} \cup \{X_k = 0 \vert \ell + 1 \leq k \leq m\}) 
\end{equation}
that is, the union of the $X_1$-,$\ldots$, $X_\ell$-coordinate axes and the affine space with as hyperplane at infinity $X = 0$.

\medskip
\subsection{The new functor, $\mF$}

The modified functor $\mF$ is defined in the same spirit as $\mS$: start with a loose graph $\Gamma$, embed it in the ambient 
projective $\Fun$-space $\mathbf{P} = \mathbf{P}(\Gamma)$, and leave out all the multiplicative groups defined by the set of edges of $\mathbf{P} \setminus \Gamma$, as 
above. More details will be provided in due course.

\begin{center}
\begin{figure}
  \begin{tikzpicture}[style=thick, scale=1.2]

\fill (0,2) circle (2pt);
\draw (-1,1) -- (0,2);
\draw (0,2) -- (1,1);

\draw (3,1) -- (4,2);
\draw (4,2) -- (5,1);

\foreach \x in {3,5}{
\fill (\x,1) circle (2pt);}
\foreach \x in {4}{
\fill (\x,2) circle (2pt);}

\draw (3,1) -- (4,2);
\draw (4,2) -- (5,1);
\draw[loosely dotted] (3,1) -- (5,1);
\end{tikzpicture}
\caption{The inclusion $A \hookrightarrow B$.}
\end{figure}
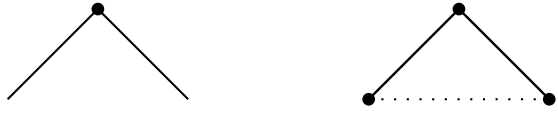
\end{center}


\section{Grothendieck ring of schemes of finte type over $\Fun$}

As we mentioned in the introduction, Deitmar schemes are based on considering commutative multiplicative monoids (with an absorbing element) as commutative rings over $\Fun$. Besides, the $\Spec$-construction allows us to have a whole scheme theory over $\Fun$ defined in an analogous way to the classical scheme theory over $\Z$. This will let us define the Grothendieck ring of schemes over $\Fun$ (Deitmar varieties).\medskip

Let us remark that the classical definition of the Grothendieck ring takes into account only varieties over a field $k$, i.e., separated integral $k$-schemes of finite type, but it is indeed possible to define the Grothendieck ring in the category of schemes of finite type over $k$ in a similar way (see \cite[section 2]{bridgeland2012}). Hence, for our study, we will directly define the Grothendieck ring for Deitmar schemes of finite type.

\begin{definition}
The Grothendieck ring of schemes of finite type over $\Fun$, denoted as $K_0(\Sch_{\Fun})$, is generated by the isomorphism classes of schemes ${X}$ of finite type over $\Fun$, $[X]_{_{\Fun}}$, with the relation
\begin{equation}
[X]_{_{\Fun}}= [X\setminus Y]_{_{\Fun}} + [Y]_{_{\Fun}} 
\end{equation}
for any closed subscheme $Y$ of $X$ and with the product structure given by
\begin{equation}
[X]_{_{\Fun}}\cdot[Y]_{_{\Fun}}= [X\times_{\Fun}Y]_{_{\Fun}}.
\end{equation}
\end{definition}
 
We denote by $\underline{\bL}=[\mathbb{A}^1_{\Fun}]_{_{\Fun}}$ the class of the affine line over $\Fun$. Let us remark that the multiplicative group $\mathbb{G}_m$ will satify, as a consequence, $[\mathbb{G}_m]_{_{\Fun}}= \underline{\bL} - 1$, since it can be identified with the affine line minus one point.\medskip
\begin{remark}
In the classical scheme theory over a field $k$, the multiplicative group $\mathbb{G}_m$ is isomorphic to the affine scheme $\Spec(k[x,y]/(xy-1))$ which has $k[x,y]/(xy - 1)$ as its coordinate ring. Nevertheless, since in the geometric theo\-ry over $\Fun$ there is only one operation on the rings, we define the kernel of a map in terms of congruences instead of ideals (see \cite[section 1.3]{lallement1979semigroups} for more details). Hence, we define the multiplicative group $\mathbb{G}_m$ over $\Fun$ to be isomorphic to the scheme $\Spec(\Fun[x, y]/(xy=1))$.
\end{remark}

\medskip
\section{Gluing}
\label{glue}

Let $\Gamma$ be a loose graph and $\mathcal{F}(\Gamma)$ be the Deitmar scheme associated to it. By definition of the functor $\mathcal{F}$ and the embedding theorem (see \S\S \ref{EmbThm}), we have for every vertex of $\Gamma$ an affine scheme of finite type over $\Fun$ defined from the loose star corresponding to the said vertex. Let us call $v_1, \ldots, v_k$ the vertices of $\Gamma$ and $\Spec(E_i)$ the affine schemes associated to $v_i$, $1\leq i\leq k$. We will study now the intersection of these affine schemes.

\begin{lemma}\label{l1}
For all $1\leq r,s, \leq k,$ $\Spec(E_r)\cap\Spec(E_s)\neq \emptyset$ if and only if $v_r$ and $v_s$ are adjacent vertices.
\end{lemma} 

\prf Suppose that $v_r$ and $v_s$ are adjacent and let $e$ be the edge having these vertices as end points. Then, $e$ belongs to the loose stars associated to $v_r$ and $v_s$ and is used to define their corresponding schemes $\Spec(E_r)$ and $\Spec(E_s)$. Thus, due to the property COV of the functor $\mathcal{F}$, the subscheme defined by $e$ is contained in the intersection.\medskip

The converse is similar. Suppose $v_r$ and $v_s$ are not adjacent. Then their respective loose stars do not have any edge in common. But the edges of the loose stars are the generators of the affine schemes $\Spec(E_r)$ and $\Spec(E_s)$, which implies that there are no relations between the generators of both schemes. Thus their intersection must be empty.\eop\medskip

This lemma implies that the gluings on the scheme $\mathcal{F}(\Gamma)$ are completely determined by the graph $\Gamma$, happening locally on the affine schemes corresponding to adjacent vertices. Besides, after Deitmar's base-extension to $\Z$, the gluings in $\mF(\Gamma)\otimes_{\Fun}\Z$ are equally determined by $\Gamma$. Let us now describe how these affine schemes are glued together.
\medskip

Consider two adjacent vertices $v_r$ and $v_s$ in $\Gamma$, with degrees $r$ and $s$ respectively, let $\Spec(E_r)$ and $\Spec(E_s)$ be their corresponding affine schemes and $e$ the edge joining both vertices.  Each of these two schemes is isomorphic to an affine space (see definition of $\mathcal{F}$, \S \ref{ppme}); then we can write $\Spec(E_r)\simeq \Spec(\Fun[x_0, \ldots, x_r])$ and $\Spec(E_s)\simeq \Spec(\Fun[y_0, \ldots, y_s])$. Since the gluing happens only at the level of the two vertices and the common edge, we can restrict ourselves to the open schemes $\mathcal{F}(\{v_r\}, \{e\})$ embedded in $\Spec(E_r)$ and $\mathcal{F}(\{v_s\}, \{e\})$ embedded in $\Spec(E_s)$, both schemes isomorphic to the affine line over $\Fun$. 
\medskip

Let us represent by $X$ the scheme $\mathcal{F}(\{v_r\}, \{e\})$ and by $Y$ the scheme $\mathcal{F}(\{v_s\}, \{e\})$. Then, w.l.o.g we can assume that $X = \Spec(\Fun[x_0])$ and $Y = \Spec(\Fun[y_0])$. We denote also by $U$ the open set $D(x_0)$ of $X$ and by $V$ the open set $D(y_0)$ of $Y$. One can easily check that the intersection $X\cap Y$ is equal to $U$ as an open set of $X$ and equal to $V$  as an open set of $Y$. Hence, to determine the gluing of $X$ and $Y$ along their intersection, we only need to define an isomorphism between the open sets $U$ and $V$.\medskip

We know that $U\simeq \Spec(\Fun[x_0, x_0^{-1}])$ and $V\simeq\Spec(\Fun[y_0, y_0^{-1}])$ and then, the isomorphism of $\Fun$-rings
\begin{center}
\begin{tabular}[column sep=2cm]{c c c}
$\Fun[x_0, x_0^{-1}]$ & $\rightarrow$ & $ \Fun[y_0, y_0^{-1}] $\\[1.2ex]
$x_0$ & $\mapsto$ & $y_0^{-1}$\\[1.2ex]
$x_0^{-1}$ & $\mapsto$ & $y_0$\\
\end{tabular}
\end{center}
induces an isomorphism of Deitmar schemes between $X$ and $Y$.\medskip

Denoting by $\iota$ the embedding of $X$ in $\Spec(E_r)$, by $\iota_1$ the embedding of $U$ in $X$ and by $j$ and $j_1$ the ones of $Y$ into $\Spec(E_s)$ and $V$ in $Y$, respectively, we can define a morphism $\psi$ of $\Fun$-schemes between $\Spec(E_r)$ and $\Spec(E_s)$ according to the following diagram
\begin{center}
\begin{tikzpicture}[auto]
\matrix(rings) [matrix of math nodes, row sep=1.2cm, column sep=2cm, ampersand replacement=\&]
{\Spec(E_r)\& \Spec(E_s) \\
\Spec(\Fun[x_0])\& \Spec(\Fun[y_0])\\
U \& V\\};
\draw[->] (rings-1-1) to node {$\psi$} (rings-1-2);
\draw[->] (rings-3-1) to node {$\simeq$}(rings-3-2);
\draw[right hook->] (rings-2-1) to node {$\iota$} (rings-1-1);
\draw[right hook->] (rings-3-1) to node {$\iota_1$} (rings-2-1);
\draw[right hook->] (rings-2-2) to node {$j$} (rings-1-2);
\draw[right hook->] (rings-3-2) to node {$j_1$} (rings-2-2);
\end{tikzpicture}
\end{center}

\noindent such that the morphism $\psi$ restricted to the intersection $\Spec(E_r)\cap\Spec(E_s)\simeq X\cap Y$ equals the isomorphism between $U$ and $V$. Doing this process for every pair of adjacent vertices in $\Gamma$, we construct morphisms of schemes compatibles with the intersections in such a way that all the affine schemes glue together to one scheme of finite type over $\Fun$.\medskip

Finally, as a direct consequence of Lemma \ref{l1}, we deduce a similar result as \cite[Theorem 4.7]{KT-Japan} for the funtor $\mathcal{F}$.

\begin{corollary}\label{conF}
Let $\Gamma$ be a loose graph and $\mathcal{F}(\Gamma)$ its loose scheme. Then, $\mathcal{F}(\Gamma)$ is connected if and only if $\Gamma$ is connected.\eop
\end{corollary}

\medskip
\section{Zeta polynomial for graphs}

Let $\Gamma$ be a loose graph. From section \S \ref{ppme}, we know that every loose graph has a Deitmar scheme $\mathcal{F}(\Gamma)$ associated to it. We will define in this section a new function on $\Gamma$ that gives us information about the class of $\mathcal{F}(\Gamma)$ in the Grothendieck ring of Deitmar schemes of finite type. To simplify, we will use the notation $[\Gamma]_{_{\Fun}}$ for the class of $\mathcal{F}(\Gamma)$ in the Grothendieck ring of schemes of finite type over $\Fun$.\medskip

We will consider first the case where $\Gamma$ is a {\em loose tree}, that is, a connected loose graph in which any two vertices are connected by an unique simple path. So, a loose tree is a connected loose graph without cycles. Note that a tree is also a loose tree.\medskip

To find this new function, we start thinking about the most basic examples, the affine spaces $\mathbb{A}_{\Fun}^n$. According to e.g. \cite{KT-Japan,Chap1}, the corresponding loose graphs are the loose graphs with one vertex (the only closed point of $\Spec(\Fun[x_1,\ldots, x_n])$) and $n$ edges having the vertex as their only point.\medskip

\begin{figure}[h]
\begin{tikzpicture}[style=thick, scale=1.2]
\draw (0,0)-- (1.5,1);
\draw (0,0)-- (1.5,0.7);
\draw (0,0)-- (1.5,-0.7);
\draw (0,0)-- (1.5,-1);
\fill (0,0) circle (2pt);
\draw[loosely dotted] (1.5, 0.7) -- (1.5, -0.7);
\draw (2, 0) node {\tiny {n-times}};
\end{tikzpicture}
\caption{Affine space $\mathbb{A}^n_{\Fun}$.}
\end{figure}
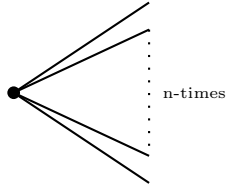

We already know that the class of the affine space $\mathbb{A}^n_{\Fun}$ in the Grothendieck ring of $\Fun$-schemes of finite type is the $n$-th power of the Lefschetz motive $\underline{\bL}$. So, our function must verify the following condition:\medskip

$\bullet$ If we denote by $\Gamma_{\mathbb{A}^n_{\Fun}}$ the loose graph corresponding to $\mathbb{A}^n_{\Fun}$, then $\big[\Gamma_{\mathbb{A}^n_{\Fun}}\big]_{_{\Fun}}= \underline{\bL}^n$.\medskip

Besides, if $\Gamma_{\mathbb{P}^1_{\Fun}}$ is the graph associated to the projective line $\mathbb{P}_{\Fun}^1$, our function should satisfy $\big[\Gamma_{\mathbb{P}^1_{\Fun}}\big]_{_{\Fun}}= \underline{\bL}+1$. 

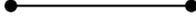
\begin{figure}[h]\label{p1}
\begin{tikzpicture}[style=thick, scale=1.2]
\draw (0,0)-- (2,0);
\fill (0,0) circle (2pt);
\fill (2,0) circle (2pt);
\end{tikzpicture}
\caption{Projective line $\mathbb{P}^1_{\Fun}$.}
\end{figure}

Comparing other examples, we came up with the following definition of our function for loose trees.

\begin{definition}
\label{D3.1}
Let $\Gamma$ be a loose tree. Then, let us consider the following notation:\medskip
\begin{itemize}
\item Let $D$ be the set of degrees $\{d_1, \ldots, d_k \}$ of $V(\Gamma)$ such that $1 < d_1 < d_2 < \ldots < d_k$.

\item Let us call $n_i$ the number of vertices of $\Gamma$ with degree $d_i$, $1\leq i \leq k$.

\item We call $\displaystyle I= \sum_{i=1}^k n_i - 1$.

\item We call $E$ the number of vertices of $\Gamma$ with degree 1, that is the {\em end points}.
\end{itemize}

Then, we define the function ``class of a loose tree,'' and we denote it as $\big[~.~\big]_{_{\Fun}}$, as follows:
\begin{center}
\begin{tabular}{c c c c}
$\big[~.~\big]_{_{\Fun}} :$ & $\{\mbox{Loose trees}\}$ & $\rightarrow$  & $K_0(\Sch_{\Fun})$\\[1.5ex]
& $\Gamma$ & $\mapsto$ & $\big[\Gamma\big]_{_{\Fun}} =  \displaystyle\sum_{i = 1}^k n_i\underline{\bL}^{d_i} - I\cdot\underline{\bL} + I + E$.
\end{tabular}
\end{center}
\end{definition}

We will prove by induction on the number of vertices of the tree that the function is well defined. We start with basic cases which will be used for proving the formula for a general loose tree.

\subsection{The projective line $\mathbb{P}^1_{\Fun}$.}

As we mentioned before, for the projective line, the corresponding tree is one edge with two end points. That gives us $E=2$, $I=-1$ and $D$ is the empty set, so the formula will be:

\begin{center}
$\big[\Gamma_{\mathbb{P}^1_{\Fun}}\big]_{_{\Fun}}= \underline{\bL}+1$,
\end{center}

\noindent which corresponds with the desired class of the projective line in the Grothendieck ring. 

\subsection{The affine space $\mathbb{A}^n_{\Fun}$.}

In this case we have one vertex of degree $n$, $E=0$ and $I=0$. Then, \begin{center}
$\big[\Gamma_{\mathbb{A}^n_{\Fun}}\big]_{_{\Fun}}=\underline{\bL}^n.$
\end{center}

\subsection{A (loose) star $S_n^k$ }

Suppose we have a star $S_n$, that means, a complete bipartite graph $K_{1,n}$.\medskip

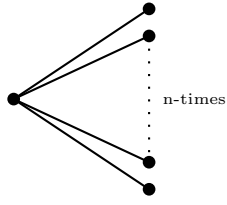
\begin{figure}[h]
\begin{tikzpicture}[style=thick, scale=1.2]
\draw (0,0)-- (1.5,1);
\draw (0,0)-- (1.5,0.7);
\draw (0,0)-- (1.5,-0.7);
\draw (0,0)-- (1.5,-1);
\fill (0,0) circle (2pt);
\fill (1.5,1) circle (2pt);
\fill (1.5,0.7) circle (2pt);
\fill (1.5,-1) circle (2pt);
\fill (1.5,-0.7) circle (2pt);
\draw[loosely dotted] (1.5, 0.7) -- (1.5, -0.7);
\draw (2, 0) node {\tiny {n-times}};
\end{tikzpicture}
\caption{A star $S_n$}
\end{figure}

We know that a vertex in the tree corresponds to a closed point of the scheme. By property ``FUCP'' (see \S\S 4.4), the disjoint union of all the vertices of degree 1 in the tree is also a closed subset of the scheme. Then, according to the relations in the Grothendieck ring, we can express the class of $S_n$ as follows.\medskip

Let us call $v_1, \ldots, v_n$ the vertices with degree one. Then
\begin{equation}
\displaystyle\big[S_n\big]_{_{\Fun}}= \big[\coprod_{i=1}^n v_i\big]_{_{\Fun}} + \big[\Gamma_{\mathbb{A}^n_{\Fun}}\big]_{_{\Fun}}= \sum_{i=1}^n \big[v_i\big]_{_{\Fun}} + \underline{\bL}^n = n\big[v_1\big]_{_{\Fun}} + \underline{\bL}^n = n + \underline{\bL}^n.
\end{equation} 

Besides, we can define a {\em loose star} with parameters $(n,k)$, $k \leq n$, and denote it as $S_n^k$, as a star $S_n$ with $k$ vertices of degree $1$. The same reasoning of above will give us a similar formula for a loose star:
\begin{equation}
\displaystyle\big[S_n^k\big]_{_{\Fun}}= k + \underline{\bL}^n.
\end{equation}

\subsection{A general loose tree}

For the case of a general loose tree we want to be able to ``break'' the tree into pieces in a certain way. To prove that the function is well defined, we will w.l.o.g. assume that a loose tree $\Gamma$ has an edge $e$ such that $\Gamma \setminus \{e\}$ is a disjoint union of two loose trees, both with at least one edge. We will call this condition $(D)$ (in case $\Gamma$ does not satisfy this condition, we will be in one of the previous cases).

\begin{lemma}
Let $\Gamma$ be a loose tree, $\mathcal{F}(\Gamma)$ the corresponding Deitmar scheme and $e$ an edge satisfying the condition $(D)$. Let $\Gamma_1$ and $\Gamma_2$ be the two disjoint loose trees obtained as above; then $\mathcal{F}(\Gamma_1)$ and $\mathcal{F}(\Gamma_2)$ are disjoint. 
\end{lemma}

\prf Since $\Gamma_1$ and $\Gamma_2$ are both connected loose trees, from Corollary 2.2, we have that $\mathcal{F}(\Gamma_1)$ and $\mathcal{F}(\Gamma_2)$ are connected Deitmar schemes. Also, by definition of the $\mathcal{F}$ functor, $\mathcal{F}(\Gamma_1 \coprod \Gamma_2)= \mathcal{F}(\Gamma_1)\cup \mathcal{F}(\Gamma_2)$ and has two connected components (Corollary 2.2). So, it is clear that $\mathcal{F}(\Gamma_1)$ and $\mathcal{F}(\Gamma_2)$ are disjoint.\eop\\

After this lemma, we are ready to prove the consistency of the formula for all trees. We will prove it using induction on trees.\medskip 

Let us start with $\Gamma$ and let $e$ be one of the vertices satisfying $(D)$. We will denote by $\bar{e}$ the subgraph having $e$ as the only edge and having two end points $v_1, v_2$. Then, $\bar{e}$ defines a projective line $\mathbb{P}^1_{\Fun}$ which is a closed subscheme of the Deitmar scheme associated to $\Gamma$. By the relations in the Grothendieck ring of schemes of finite type over $\Fun$, we have:
\begin{equation}
\big[\Gamma\big]_{_{\Fun}}= \big[\bar{e}\big]_{_{\Fun}} + \big[\Gamma\setminus\bar{e}\big]_{_{\Fun}}.
\end{equation}

\begin{remark} Let us remark here that by $[\Gamma\setminus\bar{e}]_{_{\Fun}}$ we mean the class of the Deitmar scheme defined by the loose graph $\Gamma\setminus\bar{e}~$  {\em embedded in the scheme defined by $\Gamma$}. Otherwise, if we just take the graph $\Gamma\setminus\bar{e}$ without the said embedding, we will obtain a different scheme in a projective space of higher dimension that the one in which $\mathcal{F}(\Gamma)$ is embedded!\medskip

Let's clarify this remark with the following example. Let $\Gamma$ be the graph   
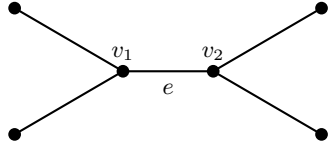
\begin{figure}[h]
\begin{tikzpicture}[style=thick, scale=1.2]
\draw (-1,0)-- (-2.2,0.7);
\draw (0,0)-- (1.2,0.7);
\draw (-1,0)-- (0,0);
\draw (0,0)-- (1.2,-0.7);
\draw (-1,0)-- (-2.2,-0.7);
\fill (0,0) circle (2pt);
\fill (-1,0) circle (2pt);
\fill (-2.2,0.7) circle (2pt);
\fill (1.2,0.7) circle (2pt);
\fill (-2.2,-0.7) circle (2pt);
\fill (1.2,-0.7) circle (2pt);
\draw (0, 0.2) node {\small {$v_2$}};
\draw (-1, 0.2) node {\small {$v_1$}};
\draw (-0.5, -0.2) node {\small {$e$}};
\end{tikzpicture}
\end{figure}

We can easily see that in this case, $\Gamma\setminus\bar{e}$ is the disjoint union of two stars $S_2$ without their respective vertices of degree $2$
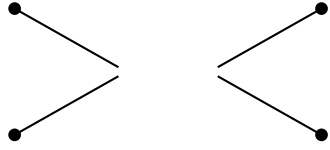
\begin{figure}[h]
\begin{tikzpicture}[style=thick, scale=1.2]
\draw (-1.05,0.05)-- (-2.2,0.7);
\draw (0.05,0.05)-- (1.2,0.7);
\draw (0.05,-0.05)-- (1.2,-0.7);
\draw (-1.05,-0.05)-- (-2.2,-0.7);
\fill (-2.2,0.7) circle (2pt);
\fill (1.2,0.7) circle (2pt);
\fill (-2.2,-0.7) circle (2pt);
\fill (1.2,-0.7) circle (2pt);
\end{tikzpicture}
\end{figure}

{\em Not} considering the scheme $\mathcal{F}(\Gamma\setminus\bar{e})$ embedded in $\mathcal{F}(\Gamma)$ gives us $4\cdot\bL$ as the zeta polynomial of $[\Gamma\setminus\bar{e}]_{_{\Fun}}$, since it corresponds to four disjoint copies of the affine line inside the projective space $\mathbb{P}_{\Fun}^7$. Nevertheless, this zeta polynomial does {\em not} give the correct information about the number of $\mathbb{F}_q$-rational points of $\mathcal{F}(\Gamma\setminus\bar{e})\otimes_{\Fun}\mathbb{F}_q$ in $\mathcal{F}(\Gamma)\otimes_{\Fun}\mathbb{F}_q$. The correct polynomial is, in fact, $2\cdot(\bL^2 + 1)$ which is the one obtained considering the right projective space in which the schemes are embedded.\end{remark}

Now, we know that $\Gamma\setminus\{e\}$ has two disjoint connected components ($\Gamma_1$ and $\Gamma_2$), and each of the end points of the edge $e$ belongs to one of these two components. Suppose that $v_1 \in \Gamma_1$ and $v_2 \in \Gamma_2$. For computing the class of $\big[\Gamma\setminus\bar{e}\big]_{_{\Fun}}$ we need to calculate the class of $\big[\Gamma_1 \setminus \{v_1\}\big]_{_{\Fun}}$ and $\big[\Gamma_2 \setminus \{v_2\}\big]_{_{\Fun}}$. To compute these two classes we also have to take into account that they define schemes embedded in $\mathcal{F}(\Gamma_1)$ and $\mathcal{F}(\Gamma_2)$, respectively.\medskip

Using the same relations on the Grothendieck ring we deduce that:
\begin{equation}
\big[\Gamma_1\cup\{e\}\cup\{v_2\}\big]_{_{\Fun}}= \big[\Gamma_1 \setminus \{v_1\}\big]_{_{\Fun}} + \big[\bar{e}\big]_{_{\Fun}}.
\end{equation}

since $\bar{e}$ defines a closed subscheme of the Deitmar scheme correponding to $\Gamma_1\cup\{e\}\cup\{v_2\}$. \medskip

Reasoning in a similar way for $\big[\Gamma_2\setminus\{v_2\}\big]_{_{\Fun}}$, we find that
\begin{equation}\label{hi}
\big[\Gamma\big]_{_{\Fun}}= \big[\bar{e}\big]_{_{\Fun}} + (\big[\Gamma_1\cup\{e\}\cup\{v_2\}\big]_{_{\Fun}} -\big[\bar{e}\big]_{_{\Fun}}) + (\big[\Gamma_2\cup\{e\}\cup\{v_1\}\big]_{_{\Fun}} -\big[\bar{e}\big]_{_{\Fun}}).
\end{equation}

Using induction, we know that the formula is valid for both $\big[\Gamma_1\cup\{e\}\cup\{v_2\}\big]_{_{\Fun}}$ and $\big[\Gamma_2\cup\{e\}\cup\{v_1\}\big]_{_{\Fun}}$. Let us write them:

\begin{equation}
\begin{aligned}
\big[\Gamma_1\cup\{e\}\cup\{v_2\}\big]_{_{\Fun}} = \displaystyle\sum_{r = 1}^{k_1} n_{1_r}\underline{\bL}^{d_{1_r}} - I_1\cdot\underline{\bL} + I_1 + E_1,\\
\big[\Gamma_2\cup\{e\}\cup\{v_1\}\big]_{_{\Fun}} = \displaystyle\sum_{j = 1}^{k_2} n_{2_j}\underline{\bL}^{d_{2_j}} - I_2\cdot\underline{\bL} + I_2 + E_2.
\end{aligned}
\end{equation}

Finally we observe that
\begin{itemize}
\item $I=I_1 + I_2 + 1$.
\item $E= E_1 +E_2 - 2$.
\item $[\bar{e}\big]_{_{\Fun}}= \underline{\bL}+1$.
\item The degree of the inner vertices and the number of vertices for each degree remain the same as in $\Gamma$.
\end{itemize}

So, introducing these two formulas in (\ref{hi}), we obtain

\begin{equation}
\displaystyle\big[\Gamma\big]_{_{\Fun}}= - (\underline{\bL} + 1) + (\sum_{r = 1}^{k_1} n_{1_r}\underline{\bL}^{d_{1_r}} + \sum_{j = 1}^{k_2} n_{2_j}\underline{\bL}^{d_{2_j}}) - (I_1 + I_2)\underline{\bL} + (I_1 + E_1 + I_2 + E_2).
\end{equation}
\medskip

Reordering the degrees of the vertices, we have
\begin{eqnarray}
\big[\Gamma\big]_{_{\Fun}} & = & \displaystyle\sum_{i = 1}^k n_i\underline{\bL}^{d_i} - (I_1 + I_2 + 1)\underline{\bL} + (I_1 + I_1 +1 + E_1 + E_2 -2)\\\
& = & \sum_{i = 1}^k n_i\underline{\bL}^{d_i} - I\cdot\underline{\bL} + I + E. \nonumber
\end{eqnarray}
\eop

\medskip
\section{Lifting $K_0(\Sch_{\Fun})$, I}
\label{lift1}

In \cite{Deitmarschemes2}, Deitmar explained how one can extend a scheme over $\Fun$ into a scheme over $\Z$ by lifting affine schemes $\Spec(A)$ to $\Spec(A)\otimes_{\Fun}\Z$, the gluing being defined by the scheme on the $\Fun$-level. The same base extension is also defined for any finite field $k$.  Thanks to the naturality of the base change functor, we will prove that this lifting is  compatible as well on the level of the Grothendieck ring of schemes of finite type.\medskip

We define $\Omega$ as a map from the Grothendieck ring $K_0(\Sch_{\Fun})$ of Deitmar schemes of finite type to the Grothendieck ring $K_0(\Sch_k)$ of schemes of finite type over any field $k$ generated by the map sending the class $\underline{\bL}$ to the class $\bL$, i.e.
\begin{eqnarray}
\Omega :  K_0(\Sch_{\Fun}) &  \rightarrow & K_0(\Sch_k)\\\
 \displaystyle\sum_{j=1}^m a_j\underline{\bL}^j & \mapsto &\sum_{j=1}^m a_j\bL^j .\nonumber
\end{eqnarray}

Let us remark that the function $\Omega$ is then defined on the subring $\Z[\underline{\bL}]$ of $K_0(\Sch_{\Fun})$.\medskip

As we did for the class of $\mathcal{F}(\Gamma)$ in the Grothendieck ring of schemes of finite type over $\Fun$, we will denote, from now on, by $[\Gamma]_{_k}$ the class of its lifting $\mathcal{F}(\Gamma)\otimes_{\Fun}k$ in the Grothendieck ring of schemes of finite type over $k$.

\begin{theorem} Let $\Gamma$ be a loose graph. Then $\Omega([\Gamma]_{_{\Fun}})=[\Gamma]_{_{k}}$.\end{theorem}

\prf We will prove this by induction, in the same way as in section 3. We will start with the basic cases.\medskip

\subsection*{Projective line} Let $\Gamma_{\mathbb{P}^1_{\Fun}}$ be the graph associated to the projective line $\mathbb{P}^1_{\Fun}$. Then, 
\begin{equation}
\Omega([\Gamma_{\mathbb{P}^1_{\Fun}}]_{_{\Fun}})= \Omega(\underline{\bL} + 1) = \bL + 1 =[\mathbb{P}_k^1]=[\Gamma_{\mathbb{P}^1_{\Fun}}]_{_k}.
\end{equation}

\subsection*{Affine spaces} Let $\Gamma_{\mathbb{A}^n_{\Fun}}$ be the graph corresponding to the affine scheme $\mathbb{A}_{\Fun}^n$. Then,
\begin{equation}
\Omega([\Gamma_{\mathbb{A}^n_{\Fun}}]_{_{\Fun}})=\Omega(\underline{\bL}^n) = \bL^n= [\mathbb{A}_{k}^n]=[\Gamma_{\mathbb{A}^n_{\Fun}}]_{_k}.
\end{equation}

\subsection*{Loose stars $S_n^k$} The Deitmar schemes associated to $S_n^k$ can be written as:
\begin{equation}
\mathcal{F}(S_n^k)= \mathbb{A}_{\Fun}^n\cup\coprod_{i=1}^k\{v_i\}.
\end{equation}

Since this $\mathcal{D}$-scheme is a disjoint union of an $n$-dimensional affine space over $\Fun$ and a disjoint union of closed points (which is closed by the property ``FUCP''), it follows that its base extension to $k$ is also a disjoint union of an $n$-dimensional $k$-affine space and a union of closed points. Hence, 
\begin{equation}
\Omega([S_n^k]_{_{\Fun}})= \Omega(\underline{\bL}^n + k) =\mathbb{L}^n + k= [\mathbb{A}^n_k]_{_k} + k= [S_n^k]_{_k}.
\end{equation}

\subsection*{Loose trees} Let $\Gamma$ be a loose tree. The way we proved the formula above for loose trees, gives us a decomposition of $[\Gamma]_{_{\Fun}}$, according to the relations in the Grothendieck ring over $\Fun$, as the sum $[\bar{e}]_{_{\Fun}} + [\Gamma \setminus \bar{e}]_{_{\Fun}}$, where $\bar{e}$ satisfies the condition $(D)$ (see \S\S 3.4).\medskip

We know that  $\bar{e}$ defines a projective line over $\Fun$, so its extension to $k$ corresponds with the projective line $\mathbb{P}^1_k$, which is also a closed subscheme of the $k$-scheme $\mathcal{F}(\Gamma)\otimes_{\Fun}k$. This gives us the following relation in the Grothendieck ring of schemes over $k$:
\begin{equation}
[\Gamma]_{_k} = [\Gamma \setminus \bar{e}]_{_k} + [\bar{e}]_{_k}.
\end{equation}

Hence, using induction on the dimesion of the loose tree and the decomposition of (\ref{hi}), we obtain that:
\begin{equation}\begin{aligned}
\Omega(\big[\Gamma_1\cup\{e\}\cup\{v_2\}\big]_{_{\Fun}})=\big[\Gamma_1\cup\{e\}\cup\{v_2\}\big]_{_k},\\
\Omega(\big[\Gamma_2\cup\{e\}\cup\{v_1\}\big]_{_{\Fun}})=\big[\Gamma_2\cup\{e\}\cup\{v_1\}\big]_{_k},\end{aligned}
\end{equation}

\noindent so, we can conclude that $\Omega([\Gamma]_{_{\Fun}})=[\Gamma]_k$.\medskip

For general loose graphs, the lifting will be obtained in \S \ref{lift2}.
\medskip

\newpage
\section{Surgery}

In this section we derive a procedure in order to inductively calculate the Grothen\-dieck polynomial of a $\mathbb{Z}$-scheme
coming from a general loose graph. In each step of the procedure, we will ``resolve'' an edge, so as to eventually end up with a tree in 
much higher dimension. So one will have to keep track of how the Grothendieck polynomial scheme changes in each step.

\medskip
\subsection{Resolution of edges}

Let $\Omega = (V,E)$ be a loose graph, and let $e \in E$ have two distinct vertices $v_1, v_2$. The {\em resolution} of $\Omega$
{\em along} $e$, denoted $\Omega_e$, is the loose graph which is obtained from $\Omega$ by deleting $e$, and adding two new
loose edges (each with one vertex) $e_1$ and $e_2$, where $v_i \in e_i$, $i = 1,2$.

One observes that 
\begin{equation}
\mathrm{dim}(\mathbf{P}(\Omega_e)) = \mathrm{dim}(\mathbf{P}(\Omega)) + 2.
\end{equation}

\medskip
\subsection{The loose graphs $\Gamma(u,v;m)$}
\label{gammauvm}

We define $\Gamma(u,v;m)$, with $m \in \mathbb{N}$ and $u, v$ symbols, to be the loose graph with adjacent
vertices $u, v$, $m$ precisely  common neighbors of $u$ and $v$ and no further incidences.

\begin{center}
\begin{figure}[h]
\label{uvm}
  \begin{tikzpicture}[style=thick, scale=1.2]
\foreach \x in {-1,1}{
\fill (\x,0) circle (2pt);}

\fill (0,1) circle (2pt);
\fill (0,1.7) circle (2pt);
\fill (0,2.4) circle (2pt);

\draw (-1,0) node[below left] {$u$} -- (1,0) node [below right] {$v$};
\draw (-1,0) -- (0,1);
\draw (0,1) -- (1,0);
\draw (-1,0) -- (0,1.7);
\draw (0,1.7) -- (1,0);
\draw (-1,0) -- (0,2.4);
\draw (0,2.4) -- (1,0);
\draw[dotted] (-1,0) -- (-0.5,1.65);
\draw[dotted] (1,0) -- (0.5,1.65);
\draw[dotted] (-1,0) -- (-0.5,1.9);
\draw[dotted] (1,0) -- (0.5,1.9);
\end{tikzpicture}
\caption{The loose graph $\Gamma(u,v;m)$.}
\end{figure}
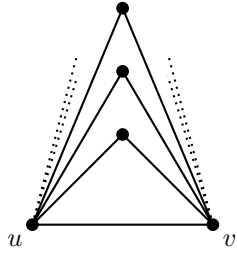
\end{center}

The corresponding $k$-schemes consist of two affine $(m + 1)$-spaces $\mathbb{A}_u$ and $\mathbb{A}_v$ and $m$ additional closed points in their spaces at infinity, of which the union covers all the points of the projective $(m + 1)$-space $\mathbf{P}(\Gamma(u,v;m))$ up to all points of the intersection $\gamma$ of their spaces at infinity (which is a projective $(m - 1)$-space), except $m$ points in $\gamma$ in general position. So the Grothendieck polynomial is 
\begin{equation}
\label{eq6}
\sum_{i = 0}^{m + 1}\uL^{i} - ((\sum_{i = 0}^{m - 1}\uL^{i}) - m) = \uL^{m + 1} + \uL^m + m.
\end{equation}

\medskip
\subsection{}
Now consider a graph $\Gamma$ which consists of a $\Gamma(u,v;m)$, $a$ further edges incident with $u$, and $b$ further edges on $v$ ($a, b \in \mathbb{N}$). We suppose that these edges are loose, but as we have see before, if they would contain some more vertices (say $c \leq a + b$ vertices), then we just add $c$ to the Grothendieck polynomial below. For further reference, denote such a loose graph by $\Gamma((u,a),(v,b);m)$.

Obviously, on the level of schemes, the only intersections occur in the projective space $\mathbf{P}((\Gamma(u,v;m)) \subseteq \mathbf{P}(\Gamma)$, so the Grothendieck polynomial is 
$$(\uL^{a + m + 1} - \uL^{m + 1}) + (\uL^{b + m + 1} - \uL^{m + 1}) + \uL^{m + 1} + \uL^m + m =$$
\begin{equation}
\uL^{a + m + 1} + \uL^{b + m + 1} - \uL^{m + 1} + \uL^m + m.
\end{equation}

\begin{remark}{\rm
Put $a = 0 = b$ and $m = 1$; then we obtain the projective $\Fun$-plane (with Grothendieck polynomial $\uL^2 + \uL  + 1$). In general, put $a = b = 0$; then we obtain (\ref{eq6}).} 
\end{remark}

\medskip
\subsection{Adding further graph structure on the common neighbors}

Define $\Gamma(u,v;G(m))$, with  $u, v$ symbols, to be the loose graph with adjacent
vertices $u, v$, $m$ common neighbors of $u$ and $v$, and with the graph $G$ defined on the common neighbors. (If the graph $G$ has no edges, we are back in \S\S \ref{gammauvm}.)  For generality's sake, let $a$ further edges be incident with $u$, and $b$ further edges with $v$ ($a, b \in \mathbb{N}$). As before, we suppose that these edges are loose.

Then in the same way as in \S\S \ref{gammauvm} one calculates
the Grothendieck polynomial to be
\begin{equation}
\uL^{a + m + 1} + \uL^{b + m + 1} - \uL^{m + 1} + \uL^m + \hP(G),
\end{equation}
where $\hP(G)$ is the Grothendieck polynomial of $G$.

\medskip
\subsection{The loose graphs $\Gamma(u,v;m)_{uv}$}

Resolving $\Gamma(u,v;m)$ along $uv$, the $k$-schemes corresponding to $\Gamma(u,v;m)_{uv}$
consist of two disjoint affine $(m + 1)$-spaces $\mathbb{A}_u$ and $\mathbb{A}_v$ (of which the hyperplanes at infinity intersect in the projective $(m - 1)$-space generated by $v_1,\ldots,v_m$) and $m$ additional mutually disjoint affine planes $\alpha_i$, $i = 1,\ldots,m$,  in the projective $(m + 3)$-space $\mathbf{P}(\Gamma(u,v;m))$ such that for each $j$, $\alpha_j \cap \mathbb{A}_u \cong \alpha_j \cap \mathbb{A}_v$ is a projective line minus two points.

The Grothendieck polynomial is 
\begin{equation}
\label{eq8}
2\uL^{m + 1} + m\uL^2 - 2m(\uL - 1).
\end{equation}

\begin{center}
\begin{figure}[h]
\label{res-uvm}
  \begin{tikzpicture}[style=thick, scale=1.2]
\foreach \x in {-1,1}{
\fill (\x,0) circle (2pt);}

\fill (0,1) circle (2pt);
\fill (0,1.7) circle (2pt);
\fill (0,2.4) circle (2pt);

\draw (-1,0) node[below left] {$u$}  -- (0,1);
\draw (0,1) -- (1,0) node [below right] {$v$};
\draw (-1,0) -- (0,1.7);
\draw (0,1.7) -- (1,0);
\draw (-1,0) -- (0,2.4);
\draw (0,2.4) -- (1,0);
\draw[dotted] (-1,0) -- (-0.5,1.65);
\draw[dotted] (1,0) -- (0.5,1.65);
\draw[dotted] (-1,0) -- (-0.5,1.9);
\draw[dotted] (1,0) -- (0.5,1.9);

\draw (-1,0) -- (-0.5,-0.5);
\draw (1,0) -- (0.5,-0.5);
\end{tikzpicture}
\caption{Resolution of $\Gamma(u,v;m)$ along the edge $uv$.}
\end{figure}
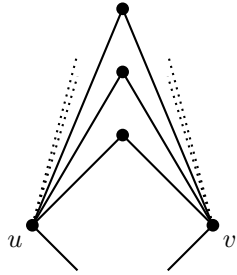
\end{center}

\medskip
\subsection{}
As before, consider a more general graph $\Gamma$, but which has $a$ loose edges on $u$, and $b$ loose edges on $v$, the case $a = 1 = b$ giving $\Gamma(u,v;m)_{uv}$. For further reference, denote such a loose graph by $\Gamma((u,a),(v,b);m)$. Then the Grothendieck polynomial is
\begin{equation}
\uL^{a + m} + \uL^{b + m} + m\uL^2 - 2m(\uL - 1).
\end{equation}

\begin{remark}{\rm
Put $a = 0 = b$ and $m = 1$; then we obtain the example studied in \S \ref{ppme} with Grothendieck polynomial $\uL^2 + 2$. If in general we put $a = b = 1$, then we obtain (\ref{eq8}).}
\end{remark}

\medskip
\subsection{Adding further graph structure on the common neighbors}

Now consider a  graph $\Gamma$ defined as in the previous subsection, but where some graph $G$ is defined on the $m$ common neighbors of $u$ and $v$.  Denote this loose graph by $\Gamma((u,a),(v,b);G)$, and remark that when $G$ has no edges, we are back in the previous subsection. 

We claim that the Grothendieck polynomial is
\begin{equation}
\uL^{a + m} + \uL^{b + m} + \hP(G)\uL^2 - 2\hP(G)(\uL - 1) = \uL^{a + m} + \uL^{b + m} + \hP(G)(\uL - 1)^2 + \hP(G).
\end{equation}

For a proof, see the next subsection.

\medskip
\subsection{General cones}

Let $\fP$ be the complete graph on $k$ vertices, and let $G_1$ and $G_2$ be subgraphs such that $\fP(G_1) \cap \fP(G_2) = \emptyset$ (noting that $\fP(G_1) \cup \fP(G_2) \subset \fP$). Define the {\em cone} with {\em base} $G_2$ and {\em vertex} $G_1$, denoted $C(G_2,G_1)$, as the  subgraph which contains $G_1$ and $G_2$, and all the edges of $\fP$ which connect a vertex of $G_2$ with a vertex of $G_1$. Note that 
$C(G_2,G_1) = C(G_1,G_2)$. 

\begin{theorem}
\label{cone}
Let the number of vertices of $G_1$ be $m_1$ and the number of vertices of $G_2$ be $m_2$. 
We have that the Grothendieck polynomial (in $K_0(\texttt{Sch}_k)$, where $k$ is any finite field, including $\Fun$) is given by
\begin{equation}
\hP(C(G_2,G_1)) = \hP(G_1)\uL^{m_2} + \hP(G_2)\uL^{m_1} - \hP(G_1)\hP(G_2)(\uL - 1).
\end{equation}
\end{theorem}

{\em Proof}.\quad
Let $k = \F_q$ be a finite field with $q$ elements, and consider $\mF(C(G_2,G_1)) \otimes_{\Fun}k$. Then for each point $\nu$
of $\mF(\hP(G_1)) \otimes_{\Fun}k$, respectively $\mF(\hP(G_2)) \otimes_{\Fun}k$, the cone contains an affine $m_2$-space defined by $\nu$ and the vertices of $\mF(\hP(G_2))\otimes_{\Fun}k$, respectively an
affine $m_1$-space defined by $\nu$ and the vertices of $\mF(\hP(G_1))\otimes_{\Fun}k$. Two by two, these affine $m_2$-spaces, respectively affine $m_1$-spaces, are disjoint, and the union of 
the affine $m_1$- and $m_2$-spaces thus obtained is the set of $k$-points of $\mF(C(G_2,G_1)) \otimes_{\Fun}k$. The order of the latter set is

\begin{equation}
\hP(G_1)(q){q}^{m_2} + \hP(G_2)(q){q}^{m_1} - \hP(G_1)(q)\hP(G_2)(q)({q} - 1),
\end{equation}
where the last term stands for the points which are doubly counted. \eop \\

Putting $G_1$ equal to a graph consisting of two vertices, one obtains the formula of the previous subsection after have added $a$ and $b$
loose edges to these vertices. 
In fact, one easily calculates the Grothendieck polynomial of a general cone with base and vertex loose graphs instead of graphs.

\begin{theorem}
\label{cone2}
Let $G_1$ and $G_2$ be loose graphs (disjoint, as above).
Let the number of vertices of $G_1$ be $m_1$ and the number of vertices of $G_2$ be $m_2$. Denote the degree of a vertex $v$ in $G_i$ by 
$\deg_{G_i}(v)$, $i = 1, 2$.
Then we have that the Grothendieck polynomial of the cone $C(G_2,G_1)$ (in $K_0(\texttt{Sch}_k)$, where $k$ is any finite field, including $\Fun$) is given by

$$\hP(C(G_2,G_1)) = \hP(G_1')\uL^{m_2} + \hP(G_2')\uL^{m_1} - \hP(G_1')\hP(G_2')(\uL - 1) $$
\begin{equation}
+ \uL^{m_2}\sum_{v \in G_1}(\uL^{\deg_{G_1}(v)} - \uL^{\deg_{G_1'}(v)}) + \uL^{m_1}\sum_{w \in G_2}(\uL^{\deg_{G_2}(w)} - \uL^{\deg_{G_2'}(w)}).
\end{equation}
\end{theorem}
{\em Proof}.\quad
Let $G_i'$ be the graph defined by $G_i$, $i = 1, 2$. By Theorem \ref{cone}, we know the Grothendieck polynomial of $C(G_2',G_1')$ in $K_0(\texttt{Sch}_k)$ for any finite field $k$. Now for each vertex $u \in G_i$, add a term $\uL^{\deg_{G_i}(u) + m_j} - \uL^{\deg_{G_i'}(u) + m_j}$,
 where $\{ i,j \} = \{ 1,2\}$. \eop \\

\medskip
\subsection{Loose graph cones versus ``classical cones'' | caution}

Let $G_1$ and $G_2$ be loose graphs, and let $k$ be a field. As above, we see $G_1$ and $G_2$ as being embedded in some 
$\Fun$-projective space, and they generate subspaces which are disjoint. 
For $A$ and $B$ two disjoint point sets in a projective space 
$\hP$ over $k$, by $A \times B$ we denote the set of points which are on lines containing a point of $A$ and a point of $B$. 

One might be tempted to think that the following identity holds:

\begin{equation}
 \mF(C(G_1,G_2)) \otimes_{\Fun}k  \cong (\mF(G_1)\otimes_{\Fun}k)\times(\mF(G_2)\otimes_{\Fun}k),
\end{equation}
that is, that the following diagram commutes:

\begin{eqnarray}
G_1, G_2\ \ \ \ \ \ &\overset{C}{\longrightarrow} &\ \ \ \ C(G_1,G_2) \nonumber \\
\ \ \nonumber \\
\ \ \downarrow{\scriptstyle (\otimes_{\Fun}k)\circ(\mF)} & &\ \ \ \ \downarrow{\scriptstyle (\otimes_{\Fun}k)\circ(\mF)} \nonumber \\  \nonumber \\
\ \ \nonumber \\
\mF(G_1)\otimes_{\Fun}k,\ \mF(G_2)\otimes_{\Fun}k  &\overset{\times}{\longrightarrow}   &(\mF(G_1)\otimes_{\Fun}k)\times(\mF(G_2)\otimes_{\Fun}k)
\nonumber \\
\end{eqnarray}

In general, this is not the case.
We will give two simple examples in which (taken that $k$ is finite), respectively,   $\vert \mF(C(G_1,G_2)) \otimes_{\Fun}k\vert > \vert (\mF(G_1)\otimes_{\Fun}k)\times(\mF(G_2)\otimes_{\Fun}k)\vert$ and  $\vert \mF(C(G_1,G_2)) \otimes_{\Fun}k\vert  < \vert (\mF(G_1)\otimes_{\Fun}k)\times(\mF(G_2)\otimes_{\Fun}k)\vert$. So there isn't even a fixed direction in which unclusion would work for 
general examples.

\medskip
\subsubsection{Example \# 1}

Let $G_2$ the graph of a projective line over $\Fun$, and let $G_1$ be the loose graph of an projective $\Fun$-plane without a multiplicative group
(i.e., an affine $\Fun$-plane with two extra points at infinity). Then the Grothendieck polynomial of $C(G_2,G_1)$ in $K_0(\texttt{Sch}_{\Fun})$ is
$\uL^4 + \uL^3 + \uL^2 + 2$, while the Grothendieck polynomial of $(\mF(G_1)\otimes_{\Fun}k)\times(\mF(G_2)\otimes_{\Fun}k)$ 
in $K_0(\texttt{Sch}_k)$ is
\begin{equation}
(\uL + 1)(\uL^2 + 2)(\uL - 1) + (\uL + 1) + (\uL^2 + 2) = \uL^4 + 2\uL^2 + \uL + 1. 
\end{equation}

Note that when the two points at infinity of $G_2$ wouldn't be there, both constructions {\em would} yield the same number of points. 
(When those two points are then added in $G_2$, the two closed points in $G_1$ suddenly see $4$-spaces instead of planes.)

\medskip
\subsubsection{Example \# 2}

Let both $G_1$ and $G_2$ be the loose graph of an affine plane over $\Fun$.
Then the Grothendieck polynomial of $C(G_2,G_1)$ (which is a tree) in $K_0(\texttt{Sch}_{\Fun})$ is
$2\uL^3 - \uL + 1$, while the Grothendieck polynomial of $(\mF(G_1)\otimes_{\Fun}k)\times(\mF(G_2)\otimes_{\Fun}k)$ 
in $K_0(\texttt{Sch}_k)$ is
\begin{equation}
\uL^2\cdot\uL^2\cdot(\uL - 1) + 2\uL^2  = \uL^5 - \uL^4 + 2\uL^2. 
\end{equation}

\medskip
\subsection{Affection principle}

Having studied a number of local situations, we now determine what happens when one resolves an edge in a general finite loose graph.
For that purpose, we consider a finite loose graph $\Gamma$, and let $\hP(\Gamma)$ be its Grothendieck polynomial in $K_0(\texttt{Sch}_{\Fun})$. We choose an edge $uv$ which is not loose, and we compare $\hP(\Gamma)$ and $\hP(\Gamma_{uv})$. 

We have seen that for each finite field $\mathbb{F}_q$, the number of $\F_q$-rational points of $\mF(\Gamma) \otimes_{\Fun} \F_q$ is given by substituting the value $q$ for the indeterminate in $\hP(\Gamma)$. As locally each closed point of $\mF(\Gamma)  \otimes_{\Fun} \F_q =: X_q$ yields an affine space (of which the dimension is the degree of the point in the graph), 
its number of points can be expressed through the Inclusion-Exclusion principle. Call the vertices of $\Gamma$ $v_1,\ldots,v_r$, and 
let for each $v_i$, $\A_i$ be the local affine space at $v_i$ of dimension $\mathrm{deg}(v_i)$. Then one can calculate the number of points (over any $\F_q$) through the expression
\begin{equation}
\sum_{i = 1}^r(-1)^{i + 1}(\sum_{1 \leq j_1 < \cdots <  j_i \leq r}\vert \A_{j_1} \cap \dots \cap \A_{j_i}\vert). 
\end{equation}

So to start with, we have to control the intersections of type $\A_x \cap \A_y$ | in other words, the intersections $\overline{\A_x} \cap
\overline{\A_y}$ (since the former intersections are controlled by the behavior at infinity). Here, $\overline{\A}$ denotes the projective completion of $\A$.\\

Calling $\rd(\cdot,\cdot)$ the distance function in $\Gamma$ defined on $V \times V$, $V$ being the vertex set (so that, for example, $\rd(s,t)$, with
$s$ and $t$ distinct vertices, is the number of edges in a shortest path from $s$ to $t$), by \S \ref{glue}, we will show we only need to consider what happens 
in the vertex set
\begin{equation}
\B(u,1) \cup \B(v,1)
\end{equation}
when resolving $uv$,
where $\B(c,k) := \{ v \in V \vert \rd(c,v) \leq k \}$. Below, we will also use the notation $v^{\perp} = \{v\}^{\perp}$ for $\{ w \in V \vert \rd(v,w) = 1\}$, and 
then $\{ v,v' \}^{\perp} := v^{\perp} \cap {v'}^{\perp}$. In particular, we have that $\B(c,1) = \{ c\} \cup c^{\perp}$.\\

The two next lemmas are immediate.

\begin{lemma}
\label{lemdist}
Let $\Gamma$ be a finite connected loose graph, and let $u$ and $v$ be distinct vertices. 
Let $k$ be any field, and consider the $k$-scheme $\mF(\Gamma) \otimes_{\Fun}k$. Then if $\rd(u,v) \geq 2$, we have that 
\begin{equation}
\A_u \cap \A_v = \emptyset.
\end{equation}
(If $\rd(u,v) \geq 3$, $\overline{\A_u} \cap \overline{\A_v} = \emptyset$.) \eop
\end{lemma}

\begin{lemma}
\label{lemdist2}
Let $\Gamma$ be a finite connected loose graph, and let $A$ be a set of distinct vertices. 
Let $k$ be any field, and consider the $k$-scheme $\mF(\Gamma) \otimes_{\Fun}k$. If $\cap_{a \in A} \A_a = \emptyset$, then 
this intersection remains empty after resolving an arbitrary edge.
\eop \\
\end{lemma}

\begin{theorem}[Affection Principle]
\label{AP}
Let $\Gamma$ be a finite connected loose graph, let $xy$ be an edge on the vertices $x$ and $y$, and let $S$ be a subset of the vertex set. 
Let $k$ be any finite field, and consider the $k$-scheme $\mF(\Gamma) \otimes_{\Fun}k$. Then $\cap_{s \in S}\A_s$ changes when one resolves 
the edge $xy$  only if $\cap_{s \in S}\A_s$ is contained in the projective subspace of $\fP(\Gamma) \otimes_{\Fun}k$ generated by $\B(x,1) \cup \B(y,1)$.
\end{theorem}
{\em Proof}.\quad
We first handle the case when $\vert S \vert = 2$, so we put $S = \{u,v\}$. As the loose edges on $u$ and $v$ clearly play no role in any change which could occur on $\A_u \cap \A_v$ when resolving $xy$, we will work w.l.o.g. in the graph $\widetilde{\Gamma}$ which is $\Gamma$ without loose edges (the ``reduced graph''); we will keep using the same notation for $x, y, u, v$. 

Clearly if $xy \not\in \fP(\B(u,1) \cup \B(v,1)) =: \fP_{u,v}$, then $\A_u \cap \A_v$ will not change through resolving $xy$, as $\overline{\A_u} \cup 
\overline{\A_v} \subseteq \fP_{u,v}$. 

Now suppose that (e.g.) $u \not\in  \B(x,1) \cup \B(y,1)$; then obviously $v \in \{ x,y\} \cup \{ x,y\}^{\perp}$ (as in each of those cases
$xy \in \overline{\A_v}$, and otherwise not). For, if $xy \in \fP_{u,v}$, then $\{x,y\} \subseteq (u^{\perp} \cup \{u\}) \cup (v^{\perp} \cup \{v\})$;
if $u \not\in  \B(x,1) \cup \B(y,1)$ ($y \not\sim u \not\sim x$ and $y \ne u \ne y$), then $x, y \in v^{\perp} \cup \{v\}$, so $v \in \{x,y\} \cup \{x,y\}^{\perp}$.

Suppose first that $v = x$. If $\rd(x,u) \geq 3$, then $\overline{\A_u} \cap \overline{\A_v} = \emptyset$ by Lemma \ref{lemdist}, and after resolution of $xy$ this stays the empty set.  If $\rd(x,u) = 2$, then $\overline{\A_u} \cap \overline{\A_v}$ remains unchanged when resolving $xy$. For, first note that $u \not\sim y$ as $u \not\in \B(x,1) \cup \B(y,1)$. Consider $\fP_{u,v} = \fP_{u,x} = \fP(\B(u,1) \cup \B(x,1))$ in $\Gamma$; as 
$\langle \overline{\A_u}, \overline{\A_x} \rangle = \fP_{u,x}$ and as $\mathrm{dim}(\fP_{u,x})$ does not change after resolving $xy$ (while in 
$\Gamma_{xy}$ we also have $\langle \overline{\A_u}, \overline{\A_x} \rangle = \fP_{u,x}$ and the dimensions of $\overline{\A_u}$, $\overline{\A_x}$ remain the same), so that  $\overline{\A_u} \cap \overline{\A_x}$ does not change ($\mathrm{dim}(\overline{\A_u} \cap \overline{\A_x})$ is in both cases determined by $\vert u^{\perp} \cap x^{\perp}\vert$).

The case $v = y$ is of course similar.

Now let $v \in x^{\perp} \cap y^{\perp}$.  By considering the cases $\rd(u,v) = 1$, $\rd(u,v) = 2$ and 
$\rd(u,v) \geq 3$, one again easily concludes that $\overline{\A_u} \cap \overline{\A_v}$ remains unchanged when resolving $xy$.
\\

We conclude that {\em if} $\overline{\A_u} \cap \overline{\A_v}$ changes after resolving $xy$, then $\{ u,v \} \subseteq  \B(x,1) \cup \B(y,1)$.\\

Now let $S \geq 3$. For obtaining a change in $\cap_{s \in S}\A_s$, at least two elements of $S$ should be in  $\B(x,1) \cup \B(y,1)$. So $\cap_{s \in S}\A_s \subseteq \fP_{x,y}$.
\eop \\

In the next corollary, $\Gamma_{\vert \fP_{u,v}}$ (e.g.) is the graph induced by $\Gamma$  in the projective space $\fP_{u,v}$ generated by $\B(u,1) \cap \B(v,1)$.

\medskip
\begin{corollary}[Geometrical Affection Principle]
\label{GAP}
Let $\Gamma$ be a finite connected loose graph, let $xy$ be an edge on the vertices $x$ and $y$, and let $k$ be any finite field.
The difference in the number of $k$-points of $\mF(\Gamma) \otimes_{\Fun}k$ and $\mF(\Gamma_{xy}) \otimes_{\Fun}k$ is 
\begin{equation}
\vert \mF(\Gamma_{\vert \fP_{x,y}}) \otimes_{\Fun}k \vert_k - \vert \mF({\Gamma_{xy}}_{\vert \fP_{x,y}}) \otimes_{\Fun}k\vert_k.
\end{equation}
In this expression, $\Gamma$ may be chosen to be reduced (but after resolving $xy$, one is of course not allowed to reduce $\Gamma_{xy}$).
\eop \\
\end{corollary}





\medskip
In terms of Grothendieck polynomials, we have the following theorem.

\begin{corollary}[Polynomial Affection principle]
Let $\Gamma$ be a finite connected loose graph, let $xy$ be an edge on the vertices $x$ and $y$, and let $k$ be any finite field.
Then in $K_0(\texttt{Sch}_k)$ we have
\begin{equation}
\hP(\Gamma) - \hP(\Gamma_{xy}) = \hP(\Gamma_{\vert \fP_{x,y}}) - \hP({\Gamma_{xy}}_{\vert \fP_{x,y}}).
\end{equation}\eop
\end{corollary}

\medskip
\subsection{Polynomial Affection principle: calculation}
\label{PAP}

Let $\Gamma$ be a finite connected loose graph, and let $e$ be an edge with vertices $x$ and $y$.
Applying the Polynomial Affection principle, we want to calculate 
\begin{equation}
\hP(\Gamma) - \hP(\Gamma_{xy}) = \hP(\Gamma_{\vert \fP_{x,y}}) - \hP({\Gamma_{xy}}_{\vert \fP_{x,y}}),
\end{equation}
in terms of certain data inside $\B(x,1) \cup \B(y,1)$.

As before, we are allowed to assume that $\Gamma$ is reduced (but not after resolution of $e$). Also, by the Affection principle, 
we only need to calculate the difference $\hP(\Gamma_{\vert \fP_{x,y}}) - \hP({\Gamma_{xy}}_{\vert \fP_{x,y}})$, so that we 
may replace $\Gamma$ by $\Gamma_{\vert \fP_{x,y}}$. We keep using the notation ``$\Gamma$'' for the latter for the sake of convenience.

\medskip
\subsection*{Before resolution}

Define $\Delta$ to be the graph (not the loose graph!) which is induced by $\Gamma$ on the vertices of $(\B(x,1) \cup \B(y,1)) \setminus \{x,y\}$, and
let $G$ be the graph (not the loose graph!) which is induced by $\Gamma$ on the vertices of $(\B(x,1) \cap \B(y,1)) \setminus \{x,y\}$. It is important to note  that $x$ and $y$ are not contained in $G$.
Let $G^L$ be the {\em loose} graph which contains  as vertex set the vertices of $G$, and with edge set the edges of $\Delta$ which contain a 
vertex of $G$ and a vertex of $\Delta$. (Note that $G$ is a subgraph of $G^L$.) 
Let $G^L_x$ be the  loose graph which has  as vertex set the vertices of $G$, and with edge set the edges of $\Delta$ which contain a 
vertex of $G$ and a vertex of $\B(x,1)$ ($\setminus \{y\}$). Define $G^L_y$ similarly, and note that $G$ is a subgraph of both 
$G^L_x$ and $G^L_y$.

\begin{center}
\begin{figure}[h]
  \begin{tikzpicture}[style=thick, scale=1.2]
\foreach \x in {-1,1}{
\fill (\x,0) circle (2pt);}

\fill (0,1) circle (2pt);
\fill (0,1.7) circle (2pt);
\fill (0,2.4) circle (2pt);

\draw (-1,1.7) circle (1.4) node [text=black,left] {$x^{\perp} \setminus \{y\}$}; 
\draw (1,1.7) circle (1.4) node [text=black,right] {$y^{\perp} \setminus \{x\}$}; 

\draw (-1,0) node[below left] {$x$} -- (1,0) node [below right] {$y$};
\draw (-1,0) -- (0,1);
\draw (0,1) -- (1,0);
\draw (-1,0) -- (0,1.7);
\draw (0,1.7) -- (1,0);
\draw (-1,0) -- (0,2.4);
\draw (0,2.4) -- (1,0);
\draw[dotted] (-1,0) -- (-0.5,1.65);
\draw[dotted] (1,0) -- (0.5,1.65);
\draw[dotted] (-1,0) -- (-0.5,1.9);
\draw[dotted] (1,0) -- (0.5,1.9);
\end{tikzpicture}
\caption{The bal $\B(x,1) \cup \B(y,1)$.}
\end{figure}
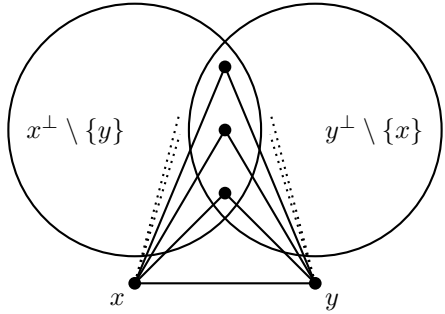
\end{center}

By the Inclusion-Exclusion principle, we know that $\hP(\Gamma)$ is given by the expression

$$\hP(\Gamma) = \hP(\A_x) + \hP(\A_y) + \hP(C(G^L,xy)) - \hP(\A_x \cap C(G^L,xy)) - \hP(\A_y \cap C(G^L,xy)) $$
\begin{equation}
\label{eqlg}
- \hP(\A_x \cap \A_y) + \hP(\A_x \cap \A_y \cap C(G^L,xy)).
\end{equation}
Here $xy$ is seen as the graph on the vertices $x$ and $y$ (i.e., a projective line).

As $\A_x \cap \A_y \subseteq C(G^L,xy)$, we have that 
\begin{equation}
\A_x \cap \A_y \cap C(G^L,xy) = \A_x \cap \A_y,
\end{equation}
so that the last two terms in Equation (\ref{eqlg}) cancel each other.

As we will see, it will be easier to calculate the following equivalent expression:

$$\hP(\Gamma) = \hP(\A_x) + \hP(\A_y) + \hP(\Delta) + (\hP(C(G^L,xy)) - \hP(G^L))$$ 
\begin{equation}
\label{eqlg2}
- \hP((\A_x \cap C(G^L,xy)) \setminus G^L) - \hP((\A_y \cap C(G^L,xy)) \setminus G^L). 
\end{equation}

\medskip
Also, we have that

\begin{equation}
\left\{ \begin{array}{ccc}
& & \\
 \hP(\A_x \cap C(G^L,xy)) \setminus G^L) &= & \hP(C(G^L_x,xy)) - \hP(C(G^L_x,y))\\
 & & \\
 \hP(\A_y \cap C(G^L,xy)) \setminus G^L) &= & \hP(C(G^L_y,xy)) - \hP(C(G^L_y,x))\\
 & & \\
 \end{array}\right..
\end{equation}




\medskip
\subsection*{After resolution}

After having resolved the edge $e$, we make a similar calculation (which is in fact a bit easier because the intersection of $\A_x$ and $\A_y$ is). We keep using the same notation.

Instead of considering the loose graph cone $C(G^L,xy)$ in the Inclusion-Exclusion principle, we only have to consider the loose graph
${G}_{(2)}^L$, which is $G^L$ with two loose edges added per vertex (one for $x$ and one for $y$), 
in order to reach all the points of $\mF(\Gamma) \otimes (\cdot)$.  So our starting point is 

$$\hP(\Gamma) = \hP(\A_x) + \hP(\A_y) + \hP(\Delta) + (\hP({G}_{(2)}^L) - \hP(G^L))$$ 
\begin{equation}
\label{eqlg3}
- \hP((\A_x \cap {G}_{(2)}^L) \setminus G^L) - \hP((\A_y \cap {G}_{(2)}^L) \setminus G^L),
\end{equation}
remarking that we immediately started with the more simple equation, and that as in the previous subsection the terms involving
$\A_x \cap \A_y$ cancel out.

\begin{center}
\begin{figure}[h]
  \begin{tikzpicture}[style=thick, scale=1.2]
\foreach \x in {-1,1}{
\fill (\x,0) circle (2pt);}

\fill (0,1) circle (2pt);
\fill (0,1.7) circle (2pt);
\fill (0,2.4) circle (2pt);

\draw (-1,1.7) circle (1.4) node [text=black,left] {$x^{\perp}$}; 
\draw (1,1.7) circle (1.4) node [text=black,right] {$y^{\perp}$}; 

\draw (-1,0) node[below left] {$x$}  -- (0,1);
\draw (0,1) -- (1,0) node [below right] {$y$};
\draw (-1,0) -- (0,1.7);
\draw (0,1.7) -- (1,0);
\draw (-1,0) -- (0,2.4);
\draw (0,2.4) -- (1,0);
\draw[dotted] (-1,0) -- (-0.5,1.65);
\draw[dotted] (1,0) -- (0.5,1.65);
\draw[dotted] (-1,0) -- (-0.5,1.9);
\draw[dotted] (1,0) -- (0.5,1.9);

\draw (-1,0) -- (-0.5,-0.5);
\draw (1,0) -- (0.5,-0.5);
\end{tikzpicture}
\caption{After resolution of the edge $xy$.}
\end{figure}
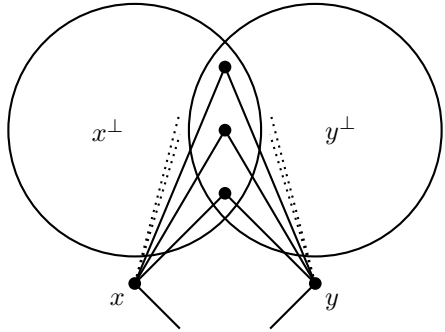
\end{center}

We immediately obtain a simple expression in the following theorem. 

\begin{theorem}[Reduction to components]
\label{compform}
Let $\{ \mC^j \vert j \in J \}$ be the set of connected components of $G^L$. Furthermore, for each $j \in J$, let $\mC^j \cap G_x^L$
be denoted by $\mC^j_x$ (and use a similar notation for $y$). Then
$$\hP(\Gamma) = \hP(\A_x) + \hP(\A_y) + \hP(\Delta) +  \sum_{j \in J}(\uL^2 - 1)\hP(\mC^j)$$
\begin{equation}
\label{eqlg4}
- \sum_{j \in J}(\uL - 1)\hP(\mC^j_x) - \sum_{j \in J}(\uL - 1)\hP(\mC^j_y).
\end{equation}
It is important to note that for each expression $\hP(\mC^\ell)$, $\hP(\mC^\ell_{x})$ and $\hP(\mC^{\ell}_y)$, one has to take 
the embedding in respectively $G^L$, $G^L_x$ and $G^L_y$ into account.
\end{theorem}

Before starting the proof, it is useful to notice that the connected components of $G^L$ correspond to the connected components of $G$. Also, note that any $\mC_x^j$ and any $\mC_y^j$ is connected.\\

{\em Proof}.\quad
First note that if $G^L$ is connected, then it follows that 
$$\hP(\Gamma) = \hP(\A_x) + \hP(\A_y) + \hP(\Delta) + (\uL^2 - 1)\hP(G^L)$$
\begin{equation}
- (\uL - 1)\hP(G^L_x) - (\uL - 1)\hP(G^L_y).
\end{equation}

Now let there be more than one connected component, and consider an arbitrary component $\mC^k$ ($k \in J$). Then $C(\{x,y\},\mC^k) \subseteq \Gamma$, and the zeta polynomial of  $C(\{x,y\},\mC^k) \setminus (\A_x \cup \A_y \cup \Delta)$ is  $(\uL^2 - 1)\hP(\mC^k)
- (\uL - 1)\hP(\mC^k_x) - (\uL - 1)\hP(\mC^k_y)$. It is easy to see that $\Gamma$ is covered by
\begin{equation}
\A_x,\ \A_y,\ \Delta,\ C(\{x,y\},\mC^j)_{\displaystyle\vert j \in J},
\end{equation}
so if we prove that for any field $k$ and $i \ne j \in J$ the following holds:
\begin{equation}
(\mF(C(\{x,y\},\mC^i))\otimes_{\Fun}k) \cap (\mF(C(\{x,y\},\mC^j))\otimes_{\Fun}k)\ \subseteq\ \A_x(k) \cup \A_y(k) \cup (\mF(\Delta)\otimes_{\Fun}k),
\end{equation}
then we are done (since each of $\A_x,\ \A_y,\ \Delta,\ C(\{x,y\},\mC^j)_{\vert j \in J}$ is contained in $\mF(\Gamma)\otimes_{\Fun}k$). Here, and below, as usual we work in the projective space $\langle \mF(\Gamma)\otimes_{\Fun}k \rangle$.

So suppose $\mC$ and $\widetilde{\mC}$ are different connected components of $G^L$, such that there is a point $z$ in 
$(\mF(C(\{x,y\},\mC))\otimes_{\Fun}k) \cap (\mF(C(\{x,y\},\widetilde{\mC}))\otimes_{\Fun}k)$ and which is not contained in $\A_x(k) \cup \A_y(k) \cup (\mF(\Delta)\otimes_{\Fun}k)$. By the structure of $\xi := \mF(C(\{x,y\},\mC))\otimes_{\Fun}k$ and $\widetilde{\xi} := \mF(C(\{x,y\},\widetilde{\mC})\otimes_{\Fun}k$, 
it follows that there are affine spaces $\mA \subseteq \xi \cap (\mF(\Delta)\otimes_{\Fun}k)$ and $\widetilde{\mA} \subseteq \widetilde{\xi}\cap (\mF(\Delta)\otimes_{\Fun}k)$, corresponding to, respectively, a vertex of $\mC$ and $\widetilde{\mC}$, such that $z \in \langle x, y, \mA\rangle \cap \langle x, y, \widetilde{\mA} \rangle$. For, $\mF(\mC)\otimes_{\Fun}k$, respectively $\mF(\widetilde{\mC})\otimes_{\Fun}k$, is covered by affine spaces $\{\A_{\mu}\}_{\vert \mu}$, respectively $\{\A_{\nu}\}_{\vert \nu}$, corresponding to vertices of $\mC$, respectively $\widetilde{\mC}$, so $\mF(C(\{x,y\},\mC))\otimes_{\Fun}k$, respectively $\mF(C(\{x,y\},\widetilde{\mC}))\otimes_{\Fun}k$, is covered by affine spaces $\{\A_{\mu, x, y}\}_{\vert \mu}$, respectively $\{\A_{\nu, x, y}\}_{\vert \nu}$, together with $\A_x$ and $\A_y$. Here $\A_{\ell, x, y}$, corresponding to the vertex $l$, denotes the affine space which contains $\A_{\ell}$ and the two extra directions defined by the edges $lx$ and $ly$. 

It follows that $z \in \mA_{x,y} \cap \widetilde{\mA}_{x,y}$ (where we use the same notation as above),  
and this is the desired contradiction, as these spaces are obviously disjoint.
\eop \\

\begin{remark}[On components]{\rm
The subdivision in connected components as in Theorem \ref{compform} is necessary. Consider for example the graph $\Gamma = (V,E)$ with $V = \{x,y,u,v,w\}$ and $E = \{ xu, xv, yu, yv, yw, wu, wv\}$. Then $G^L$ consists of the vertices $u, v$ with loose edges $uw, vw$ (embedded in a projective plane),  and hence has two connected components. If $G^L$ would be considered as being an affine plane with two additional points $u, v$ but without the vertex $w$ (that is, if $\hP(G^L) = \uL^2 + 1$), the formula
$$\hP(\Gamma) = \hP(\A_x) + \hP(\A_y) + \hP(\Delta) + (\uL^2 - 1)\hP(G^L)$$
\begin{equation}
- (\uL - 1)\hP(G^L_x) - (\uL - 1)\hP(G^L_y)
\end{equation}
would yield too many points. In other words: $G^L$ can not be considered (on the polynomial level) as being ``$\Delta$ without $w$.''
}
\end{remark}

\medskip
\subsection*{The difference}
The difference between the polynomials before and after resolution is now given by
$$\hP(\Gamma) - \hP(\Gamma_{xy}) = (\hP(C(G^L,xy)) - \hP(G^L)) 
- \hP((\A_x \cap C(G^L,xy)) \setminus G^L)$$
$$- \hP((\A_y \cap C(G^L,xy)) \setminus G^L)
- \sum_{j \in J}(\uL^2 - 1)\hP(\mC^j)$$
\begin{equation}
+ \sum_{j \in J}(\uL - 1)\hP(\mC^j_x) + \sum_{j \in J}(\uL - 1)\hP(\mC^j_y).
\end{equation}


\bigskip
\subsection{Surgery}

Recall that a {\em spanning tree} of a graph $\Omega$ is a tree which is a subgraph of $\Omega$, and has the same set of vertices as $\Omega$.
It is well known that each connected graph has at least one connected spanning tree.

\begin{theorem}
Each connected loose graph $\Gamma$ has a connected spanning tree which contains all the loose edges.
\end{theorem}
{\em Proof}.\quad
Let $S$ be the set of loose edges of $\Gamma$.
As $\Gamma$ is connected, the graph $\Gamma'$ which is obtained from $\Gamma$ by deleting all elements of $S$ is a connected
graph. So it has a connected spanning tree $T'$. Now add $S$ to $T'$ to obtain a connected spanning tree $T$ of $\Gamma$ which contains all the loose edges. \eop \\

Now let $\Gamma$ be any finite connected loose graph, and let $S$ be its set of loose edges. Then $\Gamma \setminus S =: \underline{\Gamma}$ is a finite connected graph. Let $T$ be a spanning tree of $\underline{\Gamma}$, and note that if $C$ is a cycle in $\Gamma$, it remains a cycle 
in $\underline{\Gamma}$. Define $\widetilde{S}$ to be the set of edges in $\underline{\Gamma} \setminus T$; it is by definition the set of 
{\em fundamental edges} of $\underline{\Gamma}$ with respect to the spanning tree $T$. (If one takes any edge $e \in \widetilde{S}$, then adding 
$e$ to $T$ defines a unique cycle called ``fundamental cycle.'') 

Now resolve each edge in $\widetilde{S}$ in the loose graph $\Gamma$ (in precisely $\vert \widetilde{S} \vert$ steps) as explained above, 
to obtain a tree $\overline{T}$, which contains $T$ up to a number of additional loose edges. Using the expression for the Grothedieck polynomial we calculated for trees, one now deduces the Grothendieck polynomial in $K_0(\texttt{Sch}_{\Fun})$ of $\mF(\Gamma)$.\\

\begin{theorem}
\label{Thm5.1}
Let $\Gamma$ be a loose graph and let $\underline{\Gamma}$, $T$ and $\overline{T}$ be defined as above. Then the Grothendieck polynomial in $K_0(\Sch{\Fun})$ of $\mathcal{F}(\overline{T})$ is independent of the choice of the spanning tree $T$ of $\underline{\Gamma}$.
\end{theorem}

\prf When one resolves a fundamental edge $e$ of $\underline{\Gamma}$, the edge $e$ is replaced by two different loose edges | one in each end point of $e$. Hence, the degrees of the vertices of $\Gamma$ remain the same after resolving any fundamental edge of $\underline{\Gamma}$. The result now follows from the fact that the zeta polynomial of a loose tree only depends on its spectrum of degrees (with multiplicities). \eop \\

\medskip
\medskip
\section{Lifting $K_0(\Sch_{\Fun})$, II}
\label{lift2}

We are now ready to extend the cases considered  in \S \ref{lift1} to general loose graphs.\\

First of all, we remark that we can reduce the proof of this part to the category of graphs. For, suppose we have a loose graph $\Gamma$ and let $\widetilde{\Gamma}$ be its reduced graph, i.e. the graph without loose edges. Then, it is easy to prove that the polynomials of the two graphs are the same up to a possible correction term for each vertex, a term of the form $\bL^{\mbox{\tiny deg}_{\Gamma}(v)} - \bL^{\mbox{\tiny deg}_{\widetilde{\Gamma}}(v)}$. These terms will, in the Grothendieck ring over any field $k$, locally count the points in the scheme $\mathcal{F}(\Gamma)\otimes_{\Fun} k$ that are not in the scheme $\mathcal{F}(\widetilde{\Gamma})\otimes_{\Fun} k$.\medskip

From now on, let us suppose $\Gamma$ is a graph. We will use induction on the number of vertices of the graph together with the Affection Principle and the previous cases already proved. We know that the class of the scheme defined by $\Gamma$, thanks to the Affection Principle and the surgery proccess, is indeed the class of any loose spanning tree minus the differences obtained in the resolving procedure.

\medskip
From Corollary \ref{GAP}, the difference between the zeta polynomials of the two diffe\-rent graphs (before and after resolving) is reduced to the difference between the two zeta polynomials of the subgraphs induced on the neighbors of the vertices of the resolved edge. For instance, if $xy$ is the edge on the vertices $x$ and $y$, then in $K_0(\Sch_{\Fun})$ we have
\begin{equation}\label{diff}
[\Gamma]_{\Fun} = [\Gamma_{xy}]_{\Fun} - (\underline{\bL}^2)\sum_{j \in J}[\mathcal{C}^j]_{\Fun} + (\underline{\bL} - 1)\sum_{j \in J}[\mathcal{C}^j_x]_{\Fun} + (\underline{\bL} - 1)\sum_{j \in J}[\mathcal{C}^j_y]_{\Fun} \end{equation}
\begin{equation*}
+ [C(G^L, xy)]_{\Fun}- [\A_x\cap (C(G^L, xy)\setminus G^L)]_{\Fun} - [\A_y\cap (C(G^L, xy)\setminus G^L)]_{\Fun}.
\end{equation*}\vspace{-0.2cm}

We will clarify  what the terms used in the formula mean before continuing with the proof. \medskip

The set $\{\mathcal{C}^j | j \in J\}$ is the set of connected components of $G^L$ and $\mathcal{C}^j_x$ and $\mathcal{C}^j_y$ denote $\mathcal{C}^j\cap G^L_x$ and $\mathcal{C}^j\cap G^L_y$, respectively. The following term, $C(G^L, xy)$, is the loose subgraph defined by adding to $G^L$ the vertices $x$ and $y$, the edge $xy$ and all the edges of $\Gamma$ having one end point in $\{x, y\}$ and the other in $x^{\perp}\cap y^{\perp}$. For computing the last two terms, it is sufficient to realize that, according to the relations in the Grothendieck ring, the polynomial of $\A_x\cap (C(G^L, xy)\setminus G^L)$ (respectively $\A_y\cap (C(G^L, xy)\setminus G^L)$) can be computed as the difference between the zeta polynomials of $C(G^L_x, xy)$ and $C(G^L_x, y)$ (respectively $C(G^L_y, xy)$ and $C(G^L_y, x)$). In an analogous way as before, the loose subgraph $C(G^L_x, y)$ is constructed by adding to the loose graph $G^L_x$ the vertex $y$ and all the edges of $\Gamma$ joining $y$ with the vertices of $x^{\perp}\cap y^{\perp}$. The same construction with the vertex $x$ gives the loose subgraph $C(G^L_y, x)$. \medskip

All these loose graphs come with an embedding in $\Gamma$, but due to the FUCP pro\-perty we could just as well consider their  minimal subgraphs in $\Gamma$, i.e., the subgraphs of $\Gamma$ having the same set of edges together with their corresponding embeddings. This would give indeed a suitable decomposition in the Grothendieck ring since the class of each loose subgraph is the same as the class of its corresponding minimal subgraph up to a constant regarding the number of vertices (closed points in the scheme level) added in the embedding.\medskip

We are now able to prove the theorem using a double induction on the number of vertices and the fundamental edge number. Let $\Gamma (n,t)$ be a graph (not necessarily connected) with $n$ vertices and fundamental edge number $t$. Then we claim that 
\begin{equation*}\Omega([\Gamma (n,t)]_{\Fun})=[\Gamma (n,t)]_k,\end{equation*}

\noindent for any field $k$, $n> 0$, $t\geq 0$. We take $(n,t)= (2,0)$ as the base case, since for $n=1$ the graph would be a single vertex. 

\subsubsection*{Base case {\em (2,0)}.} A graph with $2$ vertices is either a tree or two disconnected vertices. In the case of tree the theorem holds from \S \ref{lift1}. For the disconnected case the polynomial of the graph is the addition of the polynomial of the connected components. Hence, the theorem is satisfied since each component is a tree with 1 vertex.

\subsubsection*{Inductive hypothesis.} We assume the result is true for all graphs $\Gamma (r,s)$ with $r + s \leq n + t - 1$. 

\subsubsection*{Inductive step.} Let $\Gamma (n,t)$ be a graph with $n$ vertices and fundamental edge number $t$. We can assume $\Gamma(n,t)$ to be connected. Otherwise, each connected component would be a graph whose parameters satisfy the inductive hypothesis and the theorem would immediately hold. In the case of connected graphs we distinguish then three different cases: \medskip
\begin{itemize}[leftmargin=0cm, itemindent=1.25cm]
\item {\em $\Gamma (n,t)$ is a tree.} Proved by \S \ref{lift1}.
\item {\em $\Gamma (n,t)$ has a vertex $v$ of maximal degree.} Using the cone construction with $v$ as the vertex, we obtain
\begin{equation*}
\Omega([\Gamma (n,t)]_{\Fun})= \Omega(\underline{\bL}^{n-1} + [\Gamma (n,t) \setminus {v}]_{\Fun})=\bL^{n-1} + \Omega([\Gamma (n,t) \setminus {v}]_{\Fun}),
\end{equation*}

\noindent where $\Gamma (n,t)\setminus {v}$ is the subgraph of $\Gamma (n,t)$ induced on all the vertices but $v$, i.e., a graph with $n-1$ vertices that has $s \leq t$ as fundamental edge number. Hence, by the induction hypothesis, the result holds for $\Gamma (n,t)$.

\item {\em $\Gamma (n,t)$ has no vertices of maximal degree.} To prove the theorem in this case we take one fundamental edge (w.r.t. some fixed spanning tree), which we call $xy$ and we resolve it. Then, by formula (\ref{diff}), we know that:
\begin{equation*}
[\Gamma]_{\Fun} = [\Gamma_{xy}]_{\Fun} - (\underline{\bL}^2)\sum_{j \in J}[\mathcal{C}^j]_{\Fun} + (\underline{\bL} - 1)\sum_{j \in J}[\mathcal{C}^j_x]_{\Fun} + (\underline{\bL} - 1)\sum_{j \in J}[\mathcal{C}^j_y]_{\Fun} \end{equation*}\begin{equation*}
+ [C(G^L, xy)]_{\Fun} - [(C(G_x^L, xy)]_{\Fun} + [C(G_x^L, y)]_{\Fun} - [(C(G_y^L, xy)]_{\Fun} + [C(G_y^L, x)]_{\Fun}.
\end{equation*}
\end{itemize}\medskip

Thanks to the remark made above, it is sufficient to prove the theorem using the corresponding minimal graphs of the loose graphs in the formula. For all $j \in J$, the graphs $\overline{\mathcal{C}^j}$, $\overline{\mathcal{C}^j_x}$ and $\overline{\mathcal{C}^j_y}$,  have at least 2 vertices less than $\Gamma (n,t)$, so the theorem holds for them using the induction hypothesis.\medskip

For the graphs $\overline{C(G^L, xy)}$ $\overline{C(G_x^L, xy)}$, $\overline{C(G_x^L, y)}$, $\overline{C(G_y^L, xy)}$ and $\overline{C(G_y^L, x)}$ the result also holds thanks to the number of edges. Indeed, it is well know that for a graph with $n$ vertices and $e$ number of edges, the number $t$ of fundamental edges of any spanning tree equals to $e - n + 1$.  Thus, reformulating the induction hypothesis as follows\medskip

For all graphs $\Gamma$ with $r$ edges, $r \leq  e - 1$, and all $k$ field, then
\begin{equation*}\Omega([\Gamma]_{\Fun})=[\Gamma]_k, \end{equation*}

\noindent  we prove the result for all graphs $\overline{C(G^L, xy)}$ $\overline{C(G_x^L, xy)}$, $\overline{C(G_x^L, y)}$, $\overline{C(G_y^L, xy)}$ and $\overline{C(G_y^L, x)}$. By the construction of these graphs and the fact that $\Gamma(n,t)$ has no vertices of maximal degree, it is easy to check that all these graphs have less number of edges than $\Gamma(n,t)$.

There is only $\Gamma_{xy}$ remaining. Taking into account that resolving means replacing the edge $xy$ by one loose edge in each vertex $x$ and $y$, we deduce that:
\begin{equation*}
\Omega([\Gamma_{xy}]_{\Fun})= \underline{\bL}^{\mbox{\tiny deg}(x)} - \underline{\bL}^{\mbox{\tiny deg}(x) -1} + \underline{\bL}^{\mbox{\tiny deg}(y)} - \underline{\bL}^{\mbox{\tiny deg}(y)-1} + \Omega([\Gamma\setminus\{xy\}]_{\Fun}),
\end{equation*}
\noindent where $\Gamma\setminus\{xy\}$ denotes the subgraph of $\Gamma$ obtained by deleting from $\Gamma$ the edge $xy$. It is easy to check that it has the same number of vertices as $\Gamma (n,t)$ but a smaller fundamental edge number. Hence, we obtain that $\Omega([\Gamma_{xy}]_{\Fun})=[\Gamma_{xy}]_k$, for any field $k$. The theorem is now proved for any loose graph.\eop

\newpage
\section{A new zeta function for (loose) graphs}
\label{NZF}

In this final section, we are ready to introduce a new zeta function for each loose graph. 

\subsection{The zeta function}

We start with a theorem.

\begin{theorem}
For any loose graph $\Gamma$, the $\Z$-scheme $\chi := \mF(\Gamma)\otimes_{\Fun}\Z$ is defined over $\Fun$ in Kurokawa's sense. 
\end{theorem}
{\em Proof}.\quad
Let $\Gamma$ be an arbitrary finite connected loose graph. As we have seen, there exists a polynomial $P_{\chi}(X) = \sum_{i = 0}^ma_mX^m \in \Z[X]$ such that for each finite field $k = \F_q$, the number of $\F_q$-rational points is given by
\begin{equation}
N_\chi(\F_q) := \vert \mF(\Gamma) \otimes_{\Fun}\F_q \vert = P_{\chi}(q).
\end{equation}
This is precisely what we needed to prove. \eop \\

\begin{definition}[Zeta function for (loose) graphs]
Let $\Gamma$ be a loose graph, and let $\chi := \mF(\Gamma) \otimes_{\Fun}\Z$. Let $P_{\chi}(X) = \sum_{i = 0}^ma_mX^m \in \Z[X]$ be as above.
We define the {\em $\Fun$-zeta function} of $\Gamma$ as:
\begin{equation}
\zeta^{\Fun}_{\Gamma}(t) := \displaystyle \prod_{k = 0}^m(t - k)^{-a_k}.
\end{equation}
\end{definition}

In an appendix to this paper, we will compare $\zeta^{\Fun}(\cdot)$ with the Ihara zeta function for some fundamental examples of graphs.
Note that for trees, the Ihara zeta function is trivial while $\zeta^{\Fun}(\cdot)$ contains much information.

For, let $\Gamma$ be a tree. We use the notation as before:\medskip

\begin{itemize}
\item
$D$ is  the set of degrees $\{d_1, \ldots, d_m \}$ of $V(\Gamma)$ such that $1 < d_1 < d_2 < \ldots < d_m$;
\item
$n_i$ is the number of vertices of $\Gamma$ with degree $d_i$, $1\leq i \leq m$;
\item
$\displaystyle I= \sum_{i=1}^m n_i - 1$;
\item
$E$ is the number of vertices of $\Gamma$ with degree $1$.
\end{itemize}

Then, we proved that
\begin{equation}
\big[\Gamma\big]_{_{\Fun}} =  \displaystyle\sum_{i = 1}^m n_i\underline{\bL}^{d_i} - I\cdot\underline{\bL} + I + E.
\end{equation}

The zeta function is thus given by
\begin{equation}
\zeta^{\Fun}_{\Gamma}(t)\ \ =\ \ \frac{(t - 1)^I}{t^{E + I}}\cdot\displaystyle \prod_{k = 1}^m(t - k)^{-n_k}.
\end{equation}

\bigskip
\subsection{Mixed Tate motives in the Grothendieck ring}
 
One very fundamental aspect of the philosophy of $\Fun$-geometry is that the number of $\Fun$-rational points of an $\Fun$-scheme $Y$ of finite type should equal the Euler characteristic

\begin{equation}
 \chi(Y_\C) \ \ := \ \ \sum_{i=0}^{2\dim{Y}} \ (-1)^i b_i
\end{equation}
where $b_i := \dim({\mathrm{H}^i(Y_\C,\C)})$ are the Betti numbers of the complex scheme $Y_\C := Y \otimes_{\Fun}\C$. This idea is deduced from the following thought.

Let $X$ be a smooth projective scheme such that there is a polynomial $N(q) \in \Z[q]$ that counts $\F_q$-rational points, i.e.\ $\# X(\F_q)=N(q)$ for every prime power $q$. As a consequence of the comparison theorem for singular and $\ell$-adic cohomology and Deligne's proof of the Weil conjectures, we know that the counting polynomial is of the form
\begin{equation}
 N(q)\ \ = \ \ \sum_{i=0}^n b_{2i}\ q^i 
\end{equation}
and that $b_j=0$ if $j$ is odd (cf. \cite{Kurozeta}). Thus $\chi = \sum_{i=0}^n b_{2i}$ is the Euler characteristic of $X_\C$ in this case; it equals $N(1)$, which has the interpretation as the number of $\Fun$-rational points of an $\Fun$-model $X_{\Fun}$ of $X$.

Conjecturally, after the Tate conjectures, smooth projective schemes that are equipped with a counting polynomial as above, are precisely those that come with a mixed Tate motive. As our construction defines a functor 
\begin{equation}
\mF_{\Z}: \ \ \mathcal{LG} \ \ \longrightarrow \ \ \mathcal{GS}
\end{equation}
from the category of loose graphs to the category of Grothendieck schemes which all come with a counting polynomial, a natural question is whether
one can derive the Betti numbers from loose graphs $\Gamma$ for which $\mF_{\Z}$ is smooth and projective. As we will see below (in Theorem \ref{projsch}), the answer is ``yes,'' but it is also trivial, since such $\Gamma$s always give rise to projective spaces. 
 
In the next theorem, we will use the following property: 
 
\begin{itemize}
\item[INJ]
If $\Gamma$ and $\widetilde{\Gamma}$ are nonisomorphic loose graphs, then $\mF(\Gamma) \not\cong \mF(\widetilde{\Gamma})$, and for any field $k$, $\mF(\Gamma) \otimes_{\Fun}k \not\cong \mF(\widetilde{\Gamma}) \otimes_{\Fun}k$.
\end{itemize} 
 
\begin{theorem}
\label{projsch}
Let $\Gamma$ be a connected loose graph, such that $\chi = \mF(\Gamma)\otimes_{\Fun}\mathbb{Z}$ defines a (connected) projective $\mathbb{Z}$-variety. Then $\Gamma$ is a complete graph, and $\chi$ is a projective space.
\end{theorem}

{\em Proof}.\quad
Let $k$ be any field, and define $\chi_k := \mF(\Gamma)\otimes_{\Fun}k$; then $\chi_k$ is a projective variety which is embedded in a projective space $\fP$  over $k$. For any subvariety $V$ of $\fP$, define $\overline{V}$ to be its projective closure. Then from an inclusion of varieties $U \subseteq V \subseteq \fP$, we have $\overline{U} \subseteq \overline{V} \subseteq \fP$. By our construction, we know that $\chi_k$ can be covered by a set of affine spaces $\{ \mathbb{A}_v \vert v \in V(\Gamma) \}$ over $k$, all embedded in $\fP$, where $v$ runs through the vertices of $\Gamma$. As $\chi_k$ is a projective variety (and so equal to its projective closure), all points at infinity of these spaces are also contained in $\chi_k$, that is, the projective spaces they generate are also subvarieties of $\chi_k$. So $\Gamma$ is a graph.

We proceed with an induction argument on the number $\ell$ of vertices of $\Gamma$. (Obviously, the case $v = 1$ is trivial.)
Let the number of vertices of $\Gamma$ be $\ell > 1$. Then there exists a vertex, say $x$, such that the graph $\Gamma \setminus \{x\}$ (which is the graph | not the loose graph | induced on the vertex set $V(\Gamma) \setminus \{x\}$) is connected. Then for any field $k$, $\fP_x := \langle \mF(\Gamma \setminus \{x\})\otimes_{\Fun}k \rangle$ is a proper sub-projective space of $\fP$ (by COV), and as 
\begin{equation}
\mF(\Gamma \setminus \{x\})\otimes_{\Fun}k   = \fP_x \cap (\mF(\Gamma)\otimes_{\Fun}k),
\end{equation}
$\mF(\Gamma \setminus \{x\})\otimes_{\Fun}k$ is a projective variety. By induction, it is a projective space. So 
$\mF(\Gamma)\otimes_{\Fun}k$ is a union of two projective spaces, say $\fP_1$ and $\fP_2$. This means that $\Gamma$ is a (non-disjoint) union of 
two complete graphs, say $K_1$ and $K_2$ (here, we implicitly use INJ). Now let $y$ be a vertex in $K_1 \cap K_2$. Then $y$ is adjacent to all 
other vertices of $\Gamma$, so by LOC-DIM, $\mF(\Gamma)$ contains an affine $\Fun$-space of dimension $\ell - 1$. As  $\mF(\Gamma)\otimes_{\Fun}\mathbb{Z}$ defines a  projective $\mathbb{Z}$-variety, it must be a projective space of dimension $\ell - 1$.\eop \\

Note that in the previous theorem, we did not need to ask smoothness. 




\newpage
\appendix

\section{Computation of the zeta polynomial of $K_5$.}

In this appendix we will give a detailed computation of the zeta polynomial for the complete graph on 5 vertices, $K_5$. As we already know, the associated $\Fun$-scheme for $K_5$ is the projective space $\mathbb{P}^4_{\Fun}$, and so its zeta polynomial must be
\begin{equation}
[K_5] = \bL^4 + \bL^3 + \bL^2 +\bL + 1.
\end{equation}

First of all, we will choose one loose spanning tree for the graph $K_5$, i.e, a loose tree obtained after resolving all fundamental edges of a spanning tree of $K_5$. Let us remark that we can choose any loose spanning tree since its zeta polynomial is an invariant for a given graph. In our case, let us call $\Gamma$ the following loose spanning tree of $K_5$:\\
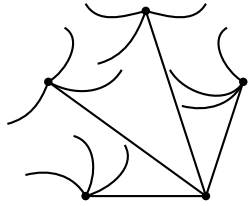
\begin{center}
\begin{tikzpicture}[style=thick, scale=0.8]
\draw (-1,0)-- (1,0);
\draw (1,0)-- (1.62,1.9);
\draw (-1.62,1.9)-- (1,0);
\draw (1,0)-- (0,3.08);
\fill (-1,0) circle (2pt);
\fill (1,0) circle (2pt);
\fill (1.62,1.9) circle (2pt);
\fill (0,3.08) circle (2pt);
\fill (-1.62,1.9) circle (2pt);
\draw (-1,0) to[out=120, in=15] (-2,0.35);
\draw (-1.62,1.9) to[out=-110, in=15] (-2.3, 1.2);
\draw (-1,0) to[out=60, in=-15] (-1.2,1);
\draw (0,3.08) to[out=-110, in=15] (-0.8,2.2);
\draw (-1,0) to[out=15, in=-60] (-0.35,0.85);
\draw (1.62,1.9) to[out=-115, in=-15] (0.6, 1.50);
\draw (0,3.08) to[out=-165, in=-60] (-1,3.2);
\draw (-1.62,1.9) to[out=30, in=-30] (-1.35, 2.8);
\draw (1.62,1.9) to[out=-120, in=-60] (0.4,2.1);
\draw (-1.62,1.9) to[out=-30, in=-120] (-0.4,2.1);
\draw (1.62,1.9) to[out=150, in=-150] (1.35,2.8);
\draw (0,3.08) to[out=-15, in=-120] (1,3.2);
\end{tikzpicture}
\end{center}

After Definition \ref{D3.1}, the zeta polynomial of $\Gamma$ can be easily computed:
\begin{equation}
[\Gamma]= 5\bL^4 - 4\bL + 4.
\end{equation}
The process to compute the zeta polynomial of the original graph $K_5$ consists of ``unresolving'' in each step one of the fundamental edges of the graph and, using the Affection Principle, keeping track of the list of resulting alternations to the zeta polynomial. By the term ``unresolving an edge of a graph $\Gamma$,'' one means choosing one graph $\Gamma_1$ in a way that $\Gamma$ is obtained from $\Gamma_1$ after resolving a fundamental edge. In our case, we take $\Gamma_1$ to be\\
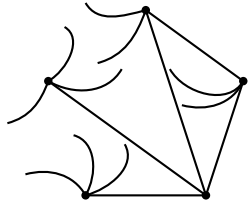
\begin{center}
\begin{tikzpicture}[style=thick, scale=0.8]
\draw (-1,0)-- (1,0);
\draw (1,0)-- (1.62,1.9);
\draw (1.62,1.9)-- (0,3.08);
\draw (-1.62,1.9)-- (1,0);
\draw (1,0)-- (0,3.08);
\fill (-1,0) circle (2pt);
\fill (1,0) circle (2pt);
\fill (1.62,1.9) circle (2pt);
\fill (0,3.08) circle (2pt);
\fill (-1.62,1.9) circle (2pt);
\draw (-1,0) to[out=120, in=15] (-2,0.35);
\draw (-1.62,1.9) to[out=-110, in=15] (-2.3, 1.2);
\draw (-1,0) to[out=60, in=-15] (-1.2,1);
\draw (0,3.08) to[out=-110, in=15] (-0.8,2.2);
\draw (-1,0) to[out=15, in=-60] (-0.35,0.85);
\draw (1.62,1.9) to[out=-115, in=-15] (0.6, 1.50);
\draw (0,3.08) to[out=-165, in=-60] (-1,3.2);
\draw (-1.62,1.9) to[out=30, in=-30] (-1.35, 2.8);
\draw (1.62,1.9) to[out=-120, in=-60] (0.4,2.1);
\draw (-1.62,1.9) to[out=-30, in=-120] (-0.4,2.1);
\end{tikzpicture}
\end{center}

Let us compare both graphs to see more clearly which are the vertices that must be taken into account according to the Affection Principle.

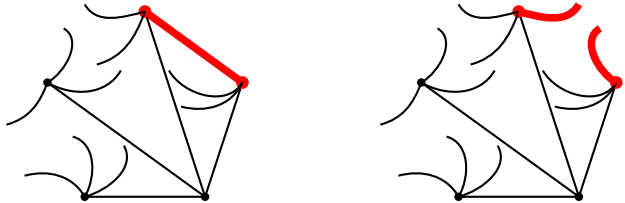
\begin{center}
\begin{tikzpicture}[style=thick, scale=0.8]
\draw (-1,0)-- (1,0);
\draw (1,0)-- (1.62,1.9);
\draw [color=red, line width=1mm](1.62,1.9) -- (0,3.08);
\draw (-1.62,1.9)-- (1,0);
\draw (1,0)-- (0,3.08);
\fill (-1,0) circle (2pt);
\fill (1,0) circle (2pt);
\fill [color=red] (1.62,1.9) circle (3pt);
\fill [color=red] (0,3.08) circle (3pt);
\fill (-1.62,1.9) circle (2pt);
\draw (-1,0) to[out=120, in=15] (-2,0.35);
\draw (-1.62,1.9) to[out=-110, in=15] (-2.3, 1.2);
\draw (-1,0) to[out=60, in=-15] (-1.2,1);
\draw (0,3.08) to[out=-110, in=15] (-0.8,2.2);
\draw (-1,0) to[out=15, in=-60] (-0.35,0.85);
\draw (1.62,1.9) to[out=-115, in=-15] (0.6, 1.50);
\draw (0,3.08) to[out=-165, in=-60] (-1,3.2);
\draw (-1.62,1.9) to[out=30, in=-30] (-1.35, 2.8);
\draw (1.62,1.9) to[out=-120, in=-60] (0.4,2.1);
\draw (-1.62,1.9) to[out=-30, in=-120] (-0.4,2.1);
\end{tikzpicture}
\hspace{1.5cm}
\begin{tikzpicture}[style=thick, scale=0.8]
\draw (-1,0)-- (1,0);
\draw (1,0)-- (1.62,1.9);
\draw (-1.62,1.9)-- (1,0);
\draw (1,0)-- (0,3.08);
\fill (-1,0) circle (2pt);
\fill (1,0) circle (2pt);
\fill [color=red](1.62,1.9) circle (3pt);
\fill [color=red](0,3.08) circle (3pt);
\fill (-1.62,1.9) circle (2pt);
\draw (-1,0) to[out=120, in=15] (-2,0.35);
\draw (-1.62,1.9) to[out=-110, in=15] (-2.3, 1.2);
\draw (-1,0) to[out=60, in=-15] (-1.2,1);
\draw (0,3.08) to[out=-110, in=15] (-0.8,2.2);
\draw (-1,0) to[out=15, in=-60] (-0.35,0.85);
\draw (1.62,1.9) to[out=-115, in=-15] (0.6, 1.50);
\draw (0,3.08) to[out=-165, in=-60] (-1,3.2);
\draw (-1.62,1.9) to[out=30, in=-30] (-1.35, 2.8);
\draw (1.62,1.9) to[out=-120, in=-60] (0.4,2.1);
\draw (-1.62,1.9) to[out=-30, in=-120] (-0.4,2.1);
\draw [color=red, line width=1mm] (1.62,1.9) to[out=150, in=-150] (1.35,2.8);
\draw [color=red, line width=1mm] (0,3.08) to[out=-15, in=-120] (1,3.2);
\end{tikzpicture}
\end{center}

The edge that has been resolved to go from $\Gamma_1$ to $\Gamma$ is the red one. Using the Affection Principle  to compute the polynomial of $\Gamma_1$, we only need to know the difference between the polynomials of the following two graphs: \\
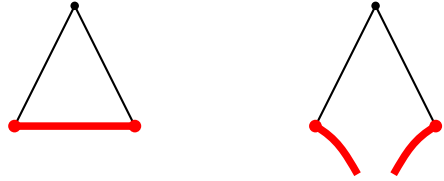
\begin{center}
\begin{tikzpicture}[style=thick, scale=0.8]
\draw [color=red, line width=1mm](-1,0)-- (1,0);
\draw (-1,0)-- (0,2);
\draw (1,0)-- (0,2);
\fill [color=red](-1,0) circle (3pt);
\fill [color=red](1,0) circle (3pt);
\fill (0,2) circle (2pt); 
\draw (4,0)-- (5,2);
\draw (6,0)-- (5,2);
\fill [color=red](4,0) circle (3pt);
\fill [color=red](6,0) circle (3pt);
\fill (5,2) circle (2pt); 
\draw [color=red, line width=1mm] (6,0) to[out=-150, in=60] (5.3,-0.8);
\draw [color=red, line width=1mm] (4,0) to[out=-30, in=120] (4.7,-0.8);
\end{tikzpicture}
\end{center}

The polynomials of these graphs are easy to compute since the graph on the left hand side is the one correponding with $\mathbb{P}^2_{\Fun}$ and the graph on the right is a tree. So the zeta polynomials are $\bL^2 + \bL + 1$ and $3\bL^2 - 2\bL + 2$, respectively. Calling $\Delta_1$ the difference between these two polynomials, we obtain the zeta polynomial for the graph $\Gamma_1$ as follows\medskip
\begin{equation}\begin{cases}
\Delta_1 = 2\bL^2 - 3\bL + 1,\\
[\Gamma_1]= [\Gamma] - \Delta_1 = 5\bL^4 - 4\bL + 4 - (2\bL^2 - 3\bL + 1)= 5\bL^4 - 2\bL^2 -\bL + 3.
\end{cases}\end{equation}\medskip

It is important to remark that, even though in this first case it was easy to compute $\Delta_1$ due to the well-known zeta polynomials of the two graphs from above, the general formulas for zeta polynomials given in \S\S \ref{PAP} would (of course) give the same expressions. In the following table, the reader can find all the necessary differences and steps which must be computed in order to get the zeta polynomial for the scheme associated to $K_5$.\medskip

\newpage
\begin{landscape}
\begin{table}[h]
\begin{tabular}{|c|c|c|c|c|c|}
\hline
 & \mbox{GRAPH}  & \mbox{BEFORE RESOLVING}  & \mbox{AFTER RESOLVING}  & $\Delta_i$ & \mbox{ZETA POLYNOMIAL} \\ \hline
$\Gamma$ & \begin{tikzpicture}[style=thick, scale=0.6]
\draw (-1,0)-- (1,0);
\draw (1,0)-- (1.62,1.9);
\draw (-1.62,1.9)-- (1,0);
\draw (1,0)-- (0,3.08);
\fill (-1,0) circle (2pt);
\fill (1,0) circle (2pt);
\fill (1.62,1.9) circle (2pt);
\fill (0,3.08) circle (2pt);
\fill (-1.62,1.9) circle (2pt);
\draw (-1,0) to[out=120, in=15] (-2,0.35);
\draw (-1.62,1.9) to[out=-110, in=15] (-2.3, 1.2);
\draw (-1,0) to[out=60, in=-15] (-1.2,1);
\draw (0,3.08) to[out=-110, in=15] (-0.8,2.2);
\draw (-1,0) to[out=15, in=-60] (-0.35,0.85);
\draw (1.62,1.9) to[out=-115, in=-15] (0.6, 1.50);
\draw (0,3.08) to[out=-165, in=-60] (-1,3.2);
\draw (-1.62,1.9) to[out=30, in=-30] (-1.35, 2.8);
\draw (1.62,1.9) to[out=-120, in=-60] (0.4,2.1);
\draw (-1.62,1.9) to[out=-30, in=-120] (-0.4,2.1);
\draw (1.62,1.9) to[out=150, in=-150] (1.35,2.8);
\draw (0,3.08) to[out=-15, in=-120] (1,3.2);
\end{tikzpicture} &  &  &  & $[\Gamma]= 5\bL^4 - 4\bL + 4$ \\ \hline
$\Gamma_1$ & \begin{tikzpicture}[style=thick, scale=0.6]
\draw (-1,0)-- (1,0);
\draw (1,0)-- (1.62,1.9);
\draw (1.62,1.9)-- (0,3.08);
\draw (-1.62,1.9)-- (1,0);
\draw (1,0)-- (0,3.08);
\fill (-1,0) circle (2pt);
\fill (1,0) circle (2pt);
\fill (1.62,1.9) circle (2pt);
\fill (0,3.08) circle (2pt);
\fill (-1.62,1.9) circle (2pt);
\draw (-1,0) to[out=120, in=15] (-2,0.35);
\draw (-1.62,1.9) to[out=-110, in=15] (-2.3, 1.2);
\draw (-1,0) to[out=60, in=-15] (-1.2,1);
\draw (0,3.08) to[out=-110, in=15] (-0.8,2.2);
\draw (-1,0) to[out=15, in=-60] (-0.35,0.85);
\draw (1.62,1.9) to[out=-115, in=-15] (0.6, 1.50);
\draw (0,3.08) to[out=-165, in=-60] (-1,3.2);
\draw (-1.62,1.9) to[out=30, in=-30] (-1.35, 2.8);
\draw (1.62,1.9) to[out=-120, in=-60] (0.4,2.1);
\draw (-1.62,1.9) to[out=-30, in=-120] (-0.4,2.1);
\end{tikzpicture} & 
\begin{tikzpicture}[style=thick, scale=0.6]
\draw (-1,0)-- (1,0);
\draw (-1,0)-- (0,2);
\draw (1,0)-- (0,2);
\fill (-1,0) circle (2pt);
\fill (1,0) circle (2pt);
\fill (0,2) circle (2pt);
\end{tikzpicture}  &
\begin{tikzpicture}[style=thick, scale=0.6]
\draw (-1,0)-- (0,2);
\draw (1,0)-- (0,2);
\fill (0,2) circle (2pt);
\fill (-1,0) circle (2pt);
\fill (1,0) circle (2pt);
\draw (1,0) to[out=-150, in=60] (0.3,-0.8);
\draw (-1,0) to[out=-30, in=120] (-0.3,-0.8);
\end{tikzpicture} &
$\Delta_1 = 2\bL^2 - 3\bL + 1$ & $[\Gamma_1]=5\bL^4 - 2\bL^2 -\bL + 3$ \\ \hline
$\Gamma_2$ & \begin{tikzpicture}[style=thick, scale=0.6]
\draw (-1,0)-- (1,0);
\draw (1,0)-- (1.62,1.9);
\draw (1.62,1.9)-- (0,3.08);
\draw (-1.62,1.9)-- (1,0);
\draw (1.62,1.9)-- (-1.62,1.9);
\draw (1,0)-- (0,3.08);
\fill (-1,0) circle (2pt);
\fill (1,0) circle (2pt);
\fill (1.62,1.9) circle (2pt);
\fill (0,3.08) circle (2pt);
\fill (-1.62,1.9) circle (2pt);
\draw (-1,0) to[out=120, in=15] (-2,0.35);
\draw (-1.62,1.9) to[out=-110, in=15] (-2.3, 1.2);
\draw (-1,0) to[out=60, in=-15] (-1.2,1);
\draw (0,3.08) to[out=-110, in=15] (-0.8,2.2);
\draw (-1,0) to[out=15, in=-60] (-0.35,0.85);
\draw (1.62,1.9) to[out=-115, in=-15] (0.6, 1.50);
\draw (0,3.08) to[out=-165, in=-60] (-1,3.2);
\draw (-1.62,1.9) to[out=30, in=-30] (-1.35, 2.8);
\end{tikzpicture} &
\begin{tikzpicture}[style=thick, scale=0.6]
\draw (-1,0)-- (1,0);
\draw (-1,0)-- (0,2);
\draw (1,0)-- (0,2);
\draw (0,2) -- (2,2);
\draw (1,0) -- (2,2);
\fill (2,2) circle (2pt);
\fill (-1,0) circle (2pt);
\fill (1,0) circle (2pt);
\fill (0,2) circle (2pt);
\end{tikzpicture}  &
\begin{tikzpicture}[style=thick, scale=0.6]
\draw (-1,0)-- (0,2);
\draw (1,0)-- (0,2);
\draw (0,2) -- (2,2);
\draw (1,0) -- (2,2);
\fill (2,2) circle (2pt);
\fill (0,2) circle (2pt);
\fill (-1,0) circle (2pt);
\fill (1,0) circle (2pt);
\draw (1,0) to[out=-150, in=60] (0.3,-0.8);
\draw (-1,0) to[out=-30, in=120] (-0.3,-0.8);
\end{tikzpicture}  &
 $\Delta_2 = \bL^3 - \bL^2$ & $[\Gamma_2] = 5\bL^4 - \bL^3 - \bL^2 -\bL +3$\\ \hline
$\Gamma_3$  & \begin{tikzpicture}[style=thick, scale=0.6]
\draw (-1.62,1.9) -- (0,3.08);
\draw (-1,0)-- (1,0);
\draw (1,0)-- (1.62,1.9);
\draw (1.62,1.9)-- (0,3.08);
\draw (-1.62,1.9)-- (1,0);
\draw (1.62,1.9)-- (-1.62,1.9);
\draw (1,0)-- (0,3.08);
\fill (-1,0) circle (2pt);
\fill (1,0) circle (2pt);
\fill (1.62,1.9) circle (2pt);
\fill (0,3.08) circle (2pt);
\fill (-1.62,1.9) circle (2pt);
\draw (-1,0) to[out=120, in=15] (-2,0.35);
\draw (-1.62,1.9) to[out=-110, in=15] (-2.3, 1.2);
\draw (-1,0) to[out=60, in=-15] (-1.2,1);
\draw (0,3.08) to[out=-110, in=15] (-0.8,2.2);
\draw (-1,0) to[out=15, in=-60] (-0.35,0.85);
\draw (1.62,1.9) to[out=-115, in=-15] (0.6, 1.50);
\end{tikzpicture} &
\begin{tikzpicture}[style=thick, scale=0.7]
\draw (-1,0)-- (1,0);
\draw (-1,0)-- (0,1);
\draw (1,0)-- (0,1);
\draw (0,1) -- (0,2);
\draw (1,0) -- (0,2);
\draw (-1,0) -- (0,2);
\fill (0,2) circle (2pt);
\fill (-1,0) circle (2pt);
\fill (1,0) circle (2pt);
\fill (0,1) circle (2pt);
\end{tikzpicture}  &
\begin{tikzpicture}[style=thick, scale=0.7]
\draw (-1,0)-- (0,1);
\draw (1,0)-- (0,1);
\draw (0,1) -- (0,2);
\draw (1,0) -- (0,2);
\draw (-1,0) -- (0,2);
\fill (0,2) circle (2pt);
\fill (-1,0) circle (2pt);
\fill (1,0) circle (2pt);
\fill (0,1) circle (2pt);
\draw (1,0) to[out=-150, in=60] (0.3,-0.8);
\draw (-1,0) to[out=-30, in=120] (-0.3,-0.8);
\end{tikzpicture}  &
$\Delta_3 = 2\bL^3 - 2\bL^2 -\bL +1$  & $[\Gamma_3] = 5\bL^4 - 3\bL^3 + \bL^2 +2$ \\ \hline
$\Gamma_4$  & \begin{tikzpicture}[style=thick, scale=0.6]
\draw (-1.62,1.9) -- (0,3.08);
\draw (-1,0)-- (1,0);
\draw (1,0)-- (1.62,1.9);
\draw (1.62,1.9)-- (0,3.08);
\draw (-1.62,1.9)-- (1,0);
\draw (-1,0)-- (1.62,1.9);
\draw (1.62,1.9)-- (-1.62,1.9);
\draw (1,0)-- (0,3.08);
\fill (-1,0) circle (2pt);
\fill (1,0) circle (2pt);
\fill (1.62,1.9) circle (2pt);
\fill (0,3.08) circle (2pt);
\fill (-1.62,1.9) circle (2pt);
\draw (-1,0) to[out=120, in=15] (-2,0.35);
\draw (-1.62,1.9) to[out=-110, in=15] (-2.3, 1.2);
\draw (-1,0) to[out=60, in=-15] (-1.2,1);
\draw (0,3.08) to[out=-110, in=15] (-0.8,2.2);
\end{tikzpicture} &
\begin{tikzpicture}[style=thick, scale=0.7]
\draw (-1,0)-- (1,0);
\draw (-1,0)-- (0,1);
\draw (1,0)-- (0,1);
\draw (0,1) -- (2,1);
\draw (1,0) -- (2,1);
\draw (1,0) -- (1,2);
\draw (0,1) -- (1,2);
\draw (2,1) -- (1,2);
\fill (1,2) circle (2pt);
\fill (2,1) circle (2pt);
\fill (-1,0) circle (2pt);
\fill (1,0) circle (2pt);
\fill (0,1) circle (2pt);
\end{tikzpicture}  &
\begin{tikzpicture}[style=thick, scale=0.7]
\draw (-1,0)-- (0,1);
\draw (1,0)-- (0,1);
\draw (0,1) -- (2,1);
\draw (1,0) -- (2,1);
\draw (1,0) -- (1,2);
\draw (0,1) -- (1,2);
\draw (2,1) -- (1,2);
\fill (1,2) circle (2pt);
\fill (2,1) circle (2pt);
\fill (-1,0) circle (2pt);
\fill (1,0) circle (2pt);
\fill (0,1) circle (2pt);
\draw (1,0) to[out=-150, in=60] (0.3,-0.8);
\draw (-1,0) to[out=-30, in=120] (-0.3,-0.8);
\end{tikzpicture} &
$\Delta_4 = \bL^4 - 2\bL^3 + 2\bL^2 -\bL$  & $[\Gamma_4] = 4\bL^4 - \bL^3 - \bL^2 +\bL +2$\\ \hline
$\Gamma_5$  & \begin{tikzpicture}[style=thick, scale=0.6]
\draw (-1.62,1.9) -- (0,3.08);
\draw (-1,0)-- (1,0);
\draw (1,0)-- (1.62,1.9);
\draw (1.62,1.9)-- (0,3.08);
\draw (0,3.08)-- (-1,0);
\draw (-1.62,1.9)-- (1,0);
\draw (-1,0)-- (1.62,1.9);
\draw (1.62,1.9)-- (-1.62,1.9);
\draw (1,0)-- (0,3.08);
\fill (-1,0) circle (2pt);
\fill (1,0) circle (2pt);
\fill (1.62,1.9) circle (2pt);
\fill (0,3.08) circle (2pt);
\fill (-1.62,1.9) circle (2pt);
\draw (-1,0) to[out=120, in=15] (-2,0.35);
\draw (-1.62,1.9) to[out=-110, in=15] (-2.3, 1.2);
\end{tikzpicture} &
\begin{tikzpicture}[style=thick, scale=0.7]
\draw (-1,0) -- (1,0);
\draw (-1,0)-- (0,1);
\draw (1,0)-- (0,1);
\draw (0,1) -- (0,2);
\draw (1,0) -- (0,2);
\draw (-1,0) -- (0,2);
\draw (2,1) -- (1,0);
\draw (2,1) -- (0,1);
\draw (2,1) -- (0,2);
\fill (2,1) circle (2pt);
\fill (0,2) circle (2pt);
\fill (-1,0) circle (2pt);
\fill (1,0) circle (2pt);
\fill (0,1) circle (2pt);
\end{tikzpicture}  &
\begin{tikzpicture}[style=thick, scale=0.7]
\draw (-1,0)-- (0,1);
\draw (1,0)-- (0,1);
\draw (0,1) -- (0,2);
\draw (1,0) -- (0,2);
\draw (-1,0) -- (0,2);
\draw (-1,0) -- (0,2);
\draw (2,1) -- (1,0);
\draw (2,1) -- (0,1);
\draw (2,1) -- (0,2);
\fill (2,1) circle (2pt);
\fill (0,2) circle (2pt);
\fill (-1,0) circle (2pt);
\fill (1,0) circle (2pt);
\fill (0,1) circle (2pt);
\draw (1,0) to[out=-150, in=60] (0.3,-0.8);
\draw (-1,0) to[out=-30, in=120] (-0.3,-0.8);
\end{tikzpicture}  &
$\Delta_5 = \bL^4 - 2\bL^2 +\bL$  & $[\Gamma_5] = 3\bL^4 - \bL^3 + \bL^2 +2$ \\ \hline
$K_5$  & \begin{tikzpicture}[style=thick, scale=0.6]
\draw (-1.62,1.9) -- (0,3.08);
\draw (-1,0)-- (1,0);
\draw (1,0)-- (1.62,1.9);
\draw (1.62,1.9)-- (0,3.08);
\draw (-1.62,1.9)-- (-1,0);
\draw (0,3.08)-- (-1,0);
\draw (-1.62,1.9)-- (1,0);
\draw (-1,0)-- (1.62,1.9);
\draw (1.62,1.9)-- (-1.62,1.9);
\draw (1,0)-- (0,3.08);
\fill (-1,0) circle (2pt);
\fill (1,0) circle (2pt);
\fill (1.62,1.9) circle (2pt);
\fill (0,3.08) circle (2pt);
\fill (-1.62,1.9) circle (2pt);
\end{tikzpicture} &
\begin{tikzpicture}[style=thick, scale=0.6]
\draw (-1.62,1.9) -- (0,3.08);
\draw (-1,0)-- (1,0);
\draw (1,0)-- (1.62,1.9);
\draw (1.62,1.9)-- (0,3.08);
\draw (-1.62,1.9)-- (-1,0);
\draw (0,3.08)-- (-1,0);
\draw (-1.62,1.9)-- (1,0);
\draw (-1,0)-- (1.62,1.9);
\draw (1.62,1.9)-- (-1.62,1.9);
\draw (1,0)-- (0,3.08);
\fill (-1,0) circle (2pt);
\fill (1,0) circle (2pt);
\fill (1.62,1.9) circle (2pt);
\fill (0,3.08) circle (2pt);
\fill (-1.62,1.9) circle (2pt);
\end{tikzpicture} &
\begin{tikzpicture}[style=thick, scale=0.6]
\draw (-1.62,1.9) -- (0,3.08);
\draw (1,0)-- (1.62,1.9);
\draw (1.62,1.9)-- (0,3.08);
\draw (-1.62,1.9)-- (-1,0);
\draw (0,3.08)-- (-1,0);
\draw (-1.62,1.9)-- (1,0);
\draw (-1,0)-- (1.62,1.9);
\draw (1.62,1.9)-- (-1.62,1.9);
\draw (1,0)-- (0,3.08);
\fill (-1,0) circle (2pt);
\fill (1,0) circle (2pt);
\fill (1.62,1.9) circle (2pt);
\fill (0,3.08) circle (2pt);
\fill (-1.62,1.9) circle (2pt);
\draw (1,0) to[out=-150, in=60] (0.3,-0.5);
\draw (-1,0) to[out=-30, in=120] (-0.3,-0.5);
\end{tikzpicture}  &
$\Delta_6 = 2\bL^4 - 2\bL^3 -\bL + 1$  & $[K_5] = \bL^4 + \bL^3 + \bL^2 + \bL +1$ \\ \hline
\end{tabular}
\end{table}

\end{landscape}

As already mentioned before, when we apply the formulas for the graphs before and after resolving, we have to take into account that the loose graphs we use are {\em embedded} in the original graph for which we are trying to obtain the zeta polynomial. We will make a detailed study of step 5 (difference between $\Gamma_4$ and $\Gamma_5$) to clarify this remark.

\section*{Computation of $[\Gamma_5]$ starting from the zeta polynomial of $\Gamma_4$.}

The polynomial associated to $\Gamma_4$ is $4\bL^4 - \bL^3 - \bL^2 +\bL +2$. Then according to the Affection Principle, in order to obtain $[\Gamma_5]$, we only have to compute the difference between the polynomials for the following two graphs:\\

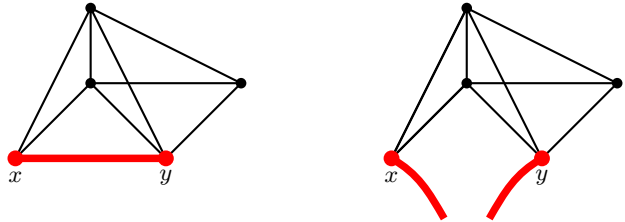
\begin{center}
\begin{tikzpicture}[style=thick, scale=1]
\draw [color=red, line width=1mm] (-1,0) -- (1,0);
\draw (-1,0)-- (0,1);
\draw (1,0)-- (0,1);
\draw (0,1) -- (0,2);
\draw (1,0) -- (0,2);
\draw (-1,0) -- (0,2);
\draw (2,1) -- (1,0);
\draw (2,1) -- (0,1);
\draw (2,1) -- (0,2);
\draw (-1, -0.25) node {\small {$x$}};
\draw (1, -0.25) node {\small {$y$}};
\fill (2,1) circle (2pt);
\fill (0,2) circle (2pt);
\fill [color=red] (-1,0) circle (3pt);
\fill [color=red] (1,0) circle (3pt);
\fill (0,1) circle (2pt);
\draw (6,0) -- (5,1);
\draw (4,0)-- (5,1);
\draw (4,0)-- (5,1);
\draw (5,1) -- (5,2);
\draw (6,0) -- (5,2);
\draw (4,0) -- (5,2);
\draw (4,0) -- (5,2);
\draw (7,1) -- (6,0);
\draw (7,1) -- (5,1);
\draw (7,1) -- (5,2);
\draw (4, -0.25) node {\small {$x$}};
\draw (6, -0.25) node {\small {$y$}};
\fill (7,1) circle (2pt);
\fill (5,2) circle (2pt);
\fill [color=red] (4,0) circle (3pt);
\fill [color=red] (6,0) circle (3pt);
\fill (5,1) circle (2pt);
\draw [color=red, line width=1mm] (6,0) to[out=-150, in=60] (5.3,-0.8);
\draw [color=red, line width=1mm] (4,0) to[out=-30, in=120] (4.7,-0.8);
\end{tikzpicture}
\end{center}

\subsection*{Before resolution.} To compute the zeta polynomial of the graph before resolution we could use the formula given in \S\S \ref{PAP} for this, but, we will use the cone construction instead | it allows us to obtain the polynomial without much computation. Thanks to the decomposition in the Grothendieck ring of varieties and the cone construction, whenever a vertex $v$ on the graph has {\em maximal} degree, i.e., equal to the number of vertices minus one, we obtain the following decomposition:\medskip
\begin{equation}
[\Gamma]=  \hP(\A_v) + [\Gamma\setminus v],
\end{equation}
\noindent where $\Gamma\setminus v$ means the ``subgraph of $\Gamma$ resulting after deleting $v$ and its incident edges.''\medskip

In our case, applying two times the cone construction, we can decompose the graph in three parts 

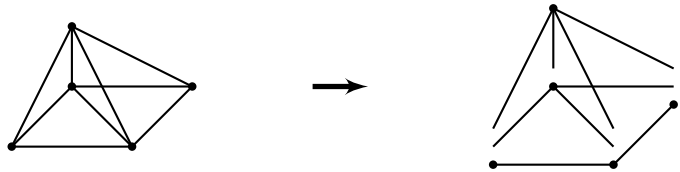
\begin{eqnarray*}
\begin{tikzpicture}[style=thick, scale=0.8]
\draw (-1,0) -- (1,0);
\draw (-1,0)-- (0,1);
\draw (1,0)-- (0,1);
\draw (0,1) -- (0,2);
\draw (1,0) -- (0,2);
\draw (-1,0) -- (0,2);
\draw (2,1) -- (1,0);
\draw (2,1) -- (0,1);
\draw (2,1) -- (0,2);
\fill (2,1) circle (2pt);
\fill (0,2) circle (2pt);
\fill (-1,0) circle (2pt);
\fill (1,0) circle (2pt);
\fill (0,1) circle (2pt);
\draw[black, -latex', line width=2pt, shorten >=1.8pt](4,1) to (5,1);
\draw (7,-0.3) -- (9,-0.3);
\draw (7,0)-- (8,1);
\draw (9,0)-- (8,1);
\draw (8,1.3) -- (8,2.3);
\draw (9,0.3) -- (8,2.3);
\draw (7,0.3) -- (8,2.3);
\draw (10,0.7) -- (9,-0.3);
\draw (10,1) -- (8,1);
\draw (10,1.3) -- (8,2.3);
\fill (10,0.7) circle (2pt);
\fill (8,2.3) circle (2pt);
\fill (7,-0.3) circle (2pt);
\fill (9,-0.3) circle (2pt);
\fill (8,1) circle (2pt);
\end{tikzpicture}
\end{eqnarray*}

\noindent and translated in terms of zeta polynomials, one obtains\medskip

\begin{equation}
[\Gamma] = \bL^4 + \bL^3 + \bL^2 + 2.
\end{equation}

\subsection*{After resolution.} According to section 7.10, the formula for the zeta polynomial in this case is\medskip
\begin{equation}\label{formula1}
\hP(\Gamma)= \hP(\A_x) + \hP(\A_y)+\hP(\Delta) + (\bL^2 - 1)\hP(G^L) - (\bL-1)\hP(G^L_x) - (\bL -1)\hP(G^L_y).
\end{equation}

Therefore, we only have to find what the different terms of the formula are ``on the graph'' Let us start with $\A_x$ and $\A_y$, which are isomorphic to $\A_{\Fun}^3$ and $\A_{\Fun}^4$, respectively. Also, observe that $\Delta$ is the subgraph generated by all the vertices except $x$ and $y$, i.e., a complete subgraph on $3$ vertices. That is, $\Delta$ defines a projective plane over $\Fun$.\medskip

Now the COV condition in the definition of the functor $\mathcal{F}$ comes into play and, as $G^L$ is embedded in $\Delta$, we deduce that it defines a projective plane minus one point; the same holds for $G^L_y$. Besides, $G^L_x$ has a projective line as associated Deitmar scheme since the only neighbors of $x$ are the common ones. Once we know all different and necessary parts of the formula (\ref{formula1}), we obtain the zeta polynomial of $\Gamma$ as follows:

\begin{align*}
\hP(\A_x)=  \bL^3, &  & \hP(\A_y)=  \bL^4, \\
\hP(\Delta)= \bL^2 + \bL + 1, &   & \hP(G^L)=   \bL^2 + \bL, \\ 
\hP(G^L_x)=  \bL +1, & & \hP(G^L_y)=  \bL^2 + \bL; 
\end{align*}

\begin{eqnarray}
\hP(\Gamma) & = & \bL^3 + \bL^4 + (\bL^2 + \bL + 1) + (\bL^2 - 1)(\bL^2 + \bL)\\ & & - (\bL - 1)(\bL^2 + \bL) - (\bL - 1)(\bL + 1)\nonumber\\
 & = & 2\bL^4 + \bL^3 - \bL^2 + \bL + 2.\nonumber
\end{eqnarray}

Hence, using the difference between the formulas obtained before and after resolution, we get that\medskip

\begin{equation*}
\left\lbrace
\begin{array}{l l l}
\Delta_5 & = &  (2\bL^4 + \bL^3 - \bL^2 + \bL + 2) - (\bL^4 + \bL^3 + \bL^2 + 2) = \bL^4 - 2\bL^2 + \bL,  \\
\hP(\Gamma_5) & = & \hP(\Gamma_4) - \Delta_5 = (4\bL^4 - \bL^3 - \bL^2 +\bL +2) - (\bL^4 - 2\bL^2 + \bL) = \\
& & \bL^4 - \bL^3 + \bL^2 + 2. 
\end{array}
\right.
\end{equation*}

\newpage
\section{Comparison with the Ihara zeta function: some examples}

In this section we will consider some examples of different graphs and compare the expressions of the corresponding Ihara zeta function and  the zeta function we have previously defined.\medskip

For the Ihara zeta function we will use the Bass-Hashimoto formula to compute its inverse (cf. \S \ref{Iharasec}) and for the new zeta function, the definition can be found in \S \ref{NZF}. In both cases we will compute the inverse so that it is easy to compare expressions. All our computations were made in Magma.\medskip

\subsection*{\textbf{Complete graph $K_4$}.}

\begin{figure}[h]\begin{tikzpicture}[style=thick, scale=1]
\draw (-1,0)-- (1,0);
\draw (-1,0)-- (0,1);
\draw (1,0)-- (0,1);
\draw (0,1) -- (0,2);
\draw (1,0) -- (0,2);
\draw (-1,0) -- (0,2);
\fill (0,2) circle (2pt);
\fill (-1,0) circle (2pt);
\fill (1,0) circle (2pt);
\fill (0,1) circle (2pt);
\end{tikzpicture}\end{figure}

A computation using the Bass-Hashimoto formula gives the following inverse of the Ihara zeta function:
\begin{equation}
\zeta (u, K_4)^{-1}= 16u^{12} - 24u^{10} -16u^{9} - 3u^{8} + 24u^{7} + 16u^{6} - 6u^{4} - 8u^{3} + 1, 
\end{equation}  

\noindent with $u \in \mathbb{C}$. For our zeta function we first need to compute the zeta polynomial, which in this case is
\begin{equation*}\mathbb{P}(K_4)= \bL^3 + \bL^2 + \bL + 1,\end{equation*}

\noindent and hence the $\Fun$-zeta function is given by
\begin{equation*}{\zeta_{K_4}^{\Fun}}^{-1}(t) = t(t-1)(t-2)(t-3) = t^4 - 6t^3 + 11t^2 - 6t.\end{equation*}

The two functions appear to be very different and this is also the case for other examples.\medskip

\subsection*{\textbf{Complete graph $K_4$ without an edge.}}Let us call  $K_4^*$ the graph $K_4 \setminus {e}$, where $e$ is any edge of $K_4$. 

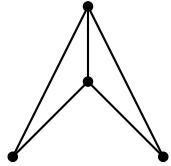
\begin{figure}[h]\begin{tikzpicture}[style=thick, scale=1]
\draw (-1,0)-- (0,1);
\draw (1,0)-- (0,1);
\draw (0,1) -- (0,2);
\draw (1,0) -- (0,2);
\draw (-1,0) -- (0,2);
\fill (0,2) circle (2pt);
\fill (-1,0) circle (2pt);
\fill (1,0) circle (2pt);
\fill (0,1) circle (2pt);
\end{tikzpicture}\end{figure}

We obtain that

\begin{eqnarray*}
\zeta (u, K_4^*)^{-1} & = & -4u^{10} + u^{8} + 4u^{7} + 4u^{6} - 2u^{4} - 4u^{3} + 1,\\
\mathbb{P}(K_4^*) & = & \bL^3 + \bL^2 + 2,\\
{\zeta_{K_4^*}^{\Fun}}^{-1}(t) & = & t^2(t-2)(t-3)= t^4 - 5t^3 + 6t^2.
\end{eqnarray*}\vspace{-0.2cm}

\subsection*{\textbf{Complete graph $K_5$.}} A detailed computation of the zeta polynomial for $K_5$ was done in the previous appendix. 

\begin{figure}[h]
\begin{tikzpicture}[style=thick, scale=1]
\draw (-1.62,1.9) -- (0,3.08);
\draw (-1,0)-- (1,0);
\draw (1,0)-- (1.62,1.9);
\draw (1.62,1.9)-- (0,3.08);
\draw (-1.62,1.9)-- (-1,0);
\draw (0,3.08)-- (-1,0);
\draw (-1.62,1.9)-- (1,0);
\draw (-1,0)-- (1.62,1.9);
\draw (1.62,1.9)-- (-1.62,1.9);
\draw (1,0)-- (0,3.08);
\fill (-1,0) circle (2pt);
\fill (1,0) circle (2pt);
\fill (1.62,1.9) circle (2pt);
\fill (0,3.08) circle (2pt);
\fill (-1.62,1.9) circle (2pt);
\end{tikzpicture}\end{figure}

In this case, the two different zeta functions are the following:

\begin{eqnarray*}
\zeta (u, K_5)^{-1} & = & -243u^{20} + 1080u^{18} + 180u^{17} - 1710u^{16} - 776u^{15} + 870u^{14} + 1200u^{13}\\
& & + 505u^{12} - 660u^{11} - 708u^{10} - 140u^{9} +165u^{8} + 240u^{7} + 70u^{6} - 24u^{5}\\
& & -30u^{4} - 20u^{3} +1,\\
\mathbb{P}(K_5) & = & \bL^4 + \bL^3 + \bL^2 + \bL + 1,\\
{\zeta_{K_5}^{\Fun}}^{-1}(t) & = & t(t-1)(t-2)(t-3)(t-4).
\end{eqnarray*}\vspace{-0.2cm}

\subsection*{\textbf{Johnson graph $J(4,2)$.}} The Johnson graph with parameters $(4,2)$ is:

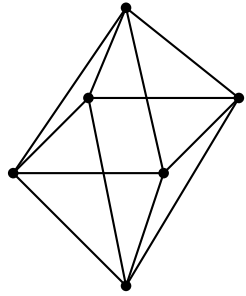
\begin{figure}[h]\begin{tikzpicture}[style=thick, scale=1]
\draw (-1,0)-- (1,0);
\draw (-1,0)-- (0,1);
\draw (0,1)-- (2,1);
\draw (2,1) -- (1,0);
\draw (-1,0) -- (0.5,2.2);
\draw (1,0) -- (0.5,2.2);
\draw (0,1) -- (0.5,2.2);
\draw (2,1) -- (0.5,2.2);
\draw (-1,0) -- (0.5,-1.5);
\draw (1,0) -- (0.5,-1.5);
\draw (0,1) -- (0.5,-1.5);
\draw (2,1) -- (0.5,-1.5);
\fill (2,1) circle (2pt);
\fill (0.5,2.2) circle (2pt);
\fill (0.5,-1.5) circle (2pt);
\fill (-1,0) circle (2pt);
\fill (1,0) circle (2pt);
\fill (0,1) circle (2pt);
\end{tikzpicture}\end{figure}

The zeta functions are given by

\begin{eqnarray*}
\zeta (u, J(4,2))^{-1} & = & 729u^{24} - 3888u^{22} - 432u^{21} + 7938u^{20} + 2160u^{19} - 6912u^{18} - 4032u^{17} \\
& & + 639u^{16} + 3008u^{15} + 2976u^{14} + 96u^{13} - 1412u^{12} - 1248u^{11} - 384u^{10} \\
& &  + 320u^9 + 327u^8 + 192u^7 + 16u^6 - 48u^5 - 30u^4 - 16u^3 + 1,\\
\mathbb{P}(J(4,2)) & = & 6\bL^4 - 12\bL^3 + 20\bL^2 -16\bL + 8,\\
{\zeta_{J(4,2)}^{\Fun}}^{-1}(t) & = & \displaystyle\frac{t^8(t-2)^{20}(t-4)^6}{(t-1)^{16}(t-3)^{12}}.
\end{eqnarray*}\vspace{-0.2cm}

\begin{remark}
Note that the Johnson graph $J(4,2)$ is the $\Fun$-analogon of the Grassmannian $\texttt{Gr}(4,2)$.
\end{remark}

\medskip
\subsection*{\textbf{Hexahedron.}} Our last example is the hexahedron, which we call $D$. Then,

\begin{figure}[h]\begin{tikzpicture}[style=thick, scale=1]
\draw (-1,0)-- (1,0);
\draw (-1,0)-- (0,1);
\draw (0,1)-- (2,1);
\draw (2,1) -- (1,0);
\draw (-1,1.5)-- (1,1.5);
\draw (-1,1.5)-- (0,2.5);
\draw (0,2.5)-- (2,2.5);
\draw (2,2.5) -- (1,1.5);
\draw (-1,0)-- (-1,1.5);
\draw (0,2.5)-- (0,1);
\draw (2,2.5)-- (2,1);
\draw (1,1.5) -- (1,0);
\fill (2,1) circle (2pt);
\fill (-1,1.5) circle (2pt);
\fill (0,2.5) circle (2pt);
\fill (2,2.5) circle (2pt);
\fill (1,1.5) circle (2pt);
\fill (-1,0) circle (2pt);
\fill (1,0) circle (2pt);
\fill (0,1) circle (2pt);
\end{tikzpicture}\end{figure}

\begin{eqnarray*}
\zeta (u, D)^{-1} & = & 256u^{24} - 768u^{22} + 480u^{20} + 400u^{18} - 183u^{16} - 384u^{14} \\
& & + 68u^{12} + 144u^{10} + 30u^8 - 32u^6 - 12u^4 + 1,\\
\mathbb{P}(D) & = & 8\bL^3 - 12\bL + 12,\\
{\zeta_{D}^{\Fun}}^{-1}(t) & = & \displaystyle\frac{t^{12}(t-3)^8}{(t-1)^{12}}.
\end{eqnarray*}

In future research we will try to see how much information about a loose graph $\Gamma$ (or its corresponding scheme $\mathcal{F}(\Gamma)\otimes_{\Fun}k$) can be obtained through $\zeta_{\{\cdot\}}^{\Fun}(t)$. We believe that the information obtained from both the new and the Ihara zeta functions could be very different.

\newpage
\section{Computations}

In this last appendix we will briefly explain the code in Magma that we use to obtain the zeta polynomial for any loose graph. We first have to remark that the main difficulty we found in writing this code is that computations are only allowed in the category of graphs. Therefore, all computations of the zeta polynomial of a loose graph $\Gamma$ will be done in what we call the {\em minimal graph} of $\Gamma$ (and denoted by $\overline{\Gamma}$), i.e. the minimal graph in which $\Gamma$ is embedded. However, computations do not get much more complicated since $\Gamma$ and $\overline{\Gamma}$ differ just by a finite number of vertices and, thanks to the FUCP property, this only implies a difference by a constant in the level of zeta polynomials. That is why, for some programs, we also need to keep track of the differences between the number of vertices.\medskip

Let us recall that to calculate the zeta polynomial of a loose graph, we use the procedure called ``surgery'' consisting of  unresolving one fundamental edge in each step and, keeping track of the differences in the graphs, ending up in a loose tree where the zeta polynomial is well defined and known. Taking all this into account, we first create a program called {\em PolTree} that given a tree computes its zeta polynomial using Definition \ref{D3.1}. However, since we need to calculate the zeta function for any loose graph and {\em PolTree} is only defined for trees, we use a program called {\em LooseSPTree} to construct a list with all different loose graphs of the surgery procedure having as its last term a loose spanning tree of the original graph.\medskip

The program {\em LoosSPTree} receives a graph $\Gamma$ as data and chooses a spanning tree of $\Gamma$ with the command {\em SpanningTree} (already implemented in Magma). Then, we compare both $\Gamma$ and its tree and we add a new graph to the output list constructed from $\Gamma$ by replacing one of the fundamental edges by two different new edges. The process continues comparing this new graph with the chosen spanning tree of $\Gamma$. Notice that in each constructed new graph we replace a fundamental edge by two loose edges. Hence, the program will finish when all the fundamental edges are resolved, i.e., when the graph constructed is a tree.\medskip

To wrap up we use the program called {\em PolSPTree} that takes a graph as input and computes the zeta polynomial of a loose spanning tree of the given graph. You can see here the code:
\begin{lstlisting}
function PolLSPTree(G)
		B:=LooseSPTree(G);
		n:=#(B);
		A:=PolTree(B[n]);
		return A;
end function;
\end{lstlisting}

Thanks to Theorem \ref{Thm5.1}, the polynomial obtained in {\em PolSPTree} is independent of the list given by {\em LooseSPTree}.\medskip

Now that we have programmed everything we need for the computation of a polynomial of a tree, we start describing how to calculate the difference in each step of the surgery. For this, let us recall the general formula for a difference after resolving the edge $xy$:

\begin{equation}\label{final}
\Delta_{xy} = (\bL^2)[G^L] - (\bL - 1)[G_x^L] - (\bL - 1)[G_y^L] - [C(G^L, xy)]\end{equation}\begin{equation*}
+ [C(G_x^L, xy)] - [C(G_x^L, y)] + [C(G_y^L, xy)] - [C(G_y^L, x)].
\end{equation*}\vspace{-0.2cm}

Thanks to the list obtained in {\em LoosSPTree} we will be able to get the new loose graphs needed to calculate the differences. For, we have created several programs called {\em GL, GLx, GLy, GL2, GL2x, GL3x, GL2y} and {\em GL3y} that compute, respectively, the loose graphs $G^L$, $G_x^L$, $G_y^L$, $C(G^L, xy)$, $C(G_x^L, xy)$, $C(G_x^L, y)$, $C(G_y^L, xy)$ and $C(G_y^L, x)$. The algorithm of all 8 programs is the same. Each of them receives two different graphs as data (the first one is the graph which arises after having resolved one edge of the second one) and, checking the neighbors of each vertex in both graphs, it identifies the resolved edge and constructs a list with the one of the forementioned loose subgraphs. There is an important fact to be remarked here regarding connectedness of the graphs. In the previous formula, the loose grahps $G^L$, $G^L_x$ and $G^L_y$ might not be connected and, in fact, the formula from above considers them as being decomposed in connected components. That is why their corresponding programs give a list where different connected components are considered as different graphs. For the other six programs this is not necessary since adding any vertex $x$ or $y$ makes the loose graphs become connected. Besides, the impossibility of working with loose graphs in Magma forces us to create new programs, {\em mL, mLx} and {\em mLy,} which will keep track, respectively, of the number of vertices of the minimal graphs $\overline{G^L}$, $\overline{G^L_x}$ and $\overline{G^L_x}$ that are not vertices of their corresponding loose graphs.\medskip

Before being able to write the final program that will compute the polynomial for any graph, we need to keep track of all different polynomials of type $G^L$, $G^L_x$, etc. in all the different steps of the surgery process. For this, we will combine both the program {\em LooseSPTree} and one program of the previous ones ({\em GL, GLx, GLy, GL2, GL2x, GL3x, GL2y} or {\em GL3y}). We only describe the program keeping track of all the $G^L$'s but there is one (exactly with the same algorithm) for each of the 8 programs listed before.  We denote by {\em ListGL} a program that takes a list of graphs, applies the program {\em GL} to any two consecutive graphs from the list and creates a new list where the elements are all the graphs obtained by the program {\em GL}.\medskip

Now we can describe the last two programs calculating the zeta polynomial for an arbitrary loose graph. The first one, {\em PolynGraph1,} receives a graph (or the minimal graph of a loose graph) and follows the following recursive algorithm. Initially we set the polynomial to be zero and we settle the two different situations in which the current program should stop: the case in which the graph is empty, giving $0$ as answer; and the case where the graph is a tree, giving {\em PolTree(G)} as the answer. This is set in the following way\medskip

\begin{lstlisting}
Q:=0;
if IsNull(G) then 
	return Q;
	else
		if IsTree(G) then
			return Q + PolTree(G);
			else
				--------------
		end if;
end if;
\end{lstlisting}			

The next step is  finding out whether the graph has a vertex of maximal degree, in which case we will use the cone construction with one vertex of maximal degree as the vertex of the cone. This is necessary to avoid the program to go in an infinite computation. For instance, if in a graph $\Gamma$, after having resolved an edge, all the other vertices are common neighbors of the two vertices from the resolved edge, then the programs {\em GL2, GL2x} and {\em GL2y} will give as result the same graph $\Gamma$ and recursion will be impossible to use. So, with the code \medskip

\begin{lstlisting}
V:=Vertices(G);
			{D:=SetToIndexedSet(Alldeg(G, #(V) -1));
			if IsEmpty(D) then
			--------------
			else 
				G1:= G - D[1];
				G2:= Components(G1);
				Q:=Q + x^(Degree(D[1])); 
				for i in [1..#(G2)] do
						Q:= Q + PolynGraph1(sub< G1 | G2[i]>);
				end for;
			end if;
\end{lstlisting}

\noindent we add, in case there is a vertex {\em v} of maximal degree, the corresponding term $\bL^{\mbox{deg}(v)}$ to the polynomial $Q$ and use recursion on the connected components of the remaining graph. Once all these cases are checked we only have to take the graph, construct the list of all intermediate steps of the surgery, compute the polynomial for the loose spanning tree (last graph of the list) and use recursion on all different graphs obtained by the programs of the type {\em ListGL} with their corresponding coefficients from the formula (\ref{final}). Let us remark that resolving edges in Magma not only means adding two edges, but also adding two more vertices, as we compute everything in the category of graphs. So when we calculate the polynomial of a spanning tree of a graph, we keep track of the number of extra vertices added in the resolving procedure. This number of extra vertices must be substracted from the polynomial of the chosen spanning tree and, in that way, it is expressed in the algortihm with the following code:

\begin{lstlisting}
L:=LooseSPTree(G);
tL:=#(L);
Q:=Q + PolLSPTree(L[tL]) - 2*(tL - 1);
\end{lstlisting}

Finally, there is only one step left to end the computation. For that purpose, we make the last program called {\em PolynGraph} that receives two arguments, a graph $\Gamma$ and a number $m$ (the number of loose edges). The program gives as output the polynomial {\em PolynGraph1$(\Gamma)$} $- ~ m$.

\begin{lstlisting}
function PolynGraph(G,m)
	Q:= PolynGraph1(G) - m;
	return Q;
end function;
\end{lstlisting}

\end{document}